\documentclass[a4paper]{article}

\usepackage{amssymb,amsthm,amsmath,amsfonts}
\usepackage{graphicx}
\usepackage{cleveref}
\usepackage{subfigure}
\numberwithin{equation}{section}

\usepackage{lscape}
\usepackage{mathrsfs}

\usepackage{algorithmic}
\usepackage[ruled]{algorithm}

\usepackage{multirow}
\usepackage{verbatim}

\usepackage{color}

\usepackage{longtable}

    \def\nb{\nonumber}
\def \Vh0{\stackrel{\circ}{V}_h}

\newcommand{\q}{\quad}

 \def\bs{\bigskip}

\usepackage{longtable}

\newtheorem{thm}{Theorem}[section]
\newtheorem{lem}[thm]{Lemma}
\newtheorem{coro}{Corollary}[section]
\newtheorem{prop}[thm]{Proposition}

\usepackage{cite}
\usepackage{enumerate}
\usepackage{url}

\title{Preconditioners and Their Analyses \\
for Edge Element Saddle-point Systems \\ Arising from
Time-harmonic Maxwell Equations}

\author{
Hua Xiang
  \thanks{ School of Mathematics and Statistics, Wuhan University, Wuhan 430072, P. R. China. The research of this author was supported by the National Natural
Science Foundation of China under grants 11571265 and 11471253.
(hxiang@whu.edu.cn).   }
\and
Shiyang Zhang
\thanks{ School of Mathematics and Statistics, Wuhan University, Wuhan 430072, P. R. China.
(hydzhang@whu.edu.cn).}
\and
  Jun Zou
   \thanks{ Department of
Mathematics,  The Chinese University of Hong Kong, Shatin, Hong Kong. The work of this author was
substantially supported by Hong Kong RGC General Research Fund
(projects 14322516 and 14306814). (zou@math.cuhk.edu.hk)}
}

\begin{document}
\date{}
\maketitle

\begin{abstract}
We shall propose and analyze some new preconditioners
for the saddle-point systems arising from the edge element discretization
of the time-harmonic Maxwell equations  in three dimensions.
We will first consider the saddle-point systems with
vanishing wave number, for which we present an important relation between the solutions of
the singular curl-curl system and the non-singular saddle-point system, then demonstrate that the saddle-point system
can be efficiently solved by the Hiptmair-Xu solver. For the saddle-point systems with non-vanishing wave numbers,
we will show that the PCG with a new preconditioner can apply for the non-singular system when wave numbers are small,
while the methods like preconditioned MINRES may apply for some existing and new preconditioners
when wave numbers are large.
The spectral behaviors of the resulting preconditioned systems for the existing and new preconditioners
are analyzed and compared, and numerical experiments are presented to demonstrate and compare
the efficiencies of  these preconditioners.
\end{abstract}


{\bf AMS subject cclassifications}: 65F10, 65N22,  65N30

\section{Introduction}\label{Introduction}
In this work we shall investigate and compare some effective preconditioning solvers for
the following saddle-point system:
\begin{equation}
\label{eq.7}
 \mathcal{K}  \left(\begin{array}{c}u\\ p\end{array}\right)
 \equiv
 \begin{pmatrix}
A-k^2M & B^T \\
B & 0
\end{pmatrix} \left(\begin{array}{c}u\\ p\end{array}\right)=\left(\begin{array}{c}f\\ g\end{array}\right)
\end{equation}
where $u\in \mathbb{R} ^{n}$, $p\in \mathbb{R} ^{m}$,
$A, M\in \mathbb{R} ^{n \times n} $ and $B\in \mathbb{R} ^{m \times n}$, with $m \leq n$.
We assume that ${\cal K}$ is nonsingular, so $B$ must be
of full row rank. We are particularly interested in the case where $A$ is symmetric semi-positive definite,
and dim$(ker(A))=m$, that is,  $A$ is maximally rank deficient  \cite{Greif2015} \cite{Greif2007}.
The matrix $M$ is assumed to be symmetric positive definite, and $k$ is a given real number.
%

The saddle-point system of form \eqref{eq.7} with a maximal rank deficient $A$
arises from many applications, including the numerical solution
of time-harmonic Maxwell's equations \cite{Greif2007, Hiptmair2002, Perugia2002a}
where $k$ represents the wave number,
the underdetermined norm-minimization problems and geophysical inverse problems; see more details
in the recent paper \cite{Greif2015} that was a very inspiring and innovative work and developed
a class of indefinite block preconditioners for the use with the CG method,
and CG may converge rapidly under certain conditions when it is applied for solving the
general saddle-point system of form \eqref{eq.7} with a maximal rank deficient $A$ and $k=0$.
It was also pointed out in \cite{Greif2015} that the saddle-point system under the above
particular setting has not received as much attention as other situations, for
example, the case of a symmetric positive definite $A$. But as it was demonstrated in \cite{Greif2015}
for the special case $k=0$,
when $A$ is maximally rank deficient, some nice mathematical structures may be revealed and adopted
to help construct efficient solution methods. This work is initiated and motivated by \cite{Greif2015}
and will develop further in this direction, and show that
new efficient numerical methods can be equally constructed for more
general and difficult case $k\ne 0$.

Though most results of this work apply also to the general saddle-point system of
form \eqref{eq.7} with a maximal rank deficient $A$ as it was done in \cite{Greif2015}
for the case $k=0$, we shall mainly focus on the saddle-point system \eqref{eq.7}
that arises from the edge element discretization of the following time-harmonic Maxwell equations
in vacuum \cite{Perugia2002a,Houston2005,Chen:2000:FEM:588310.588552,Demkowicz1998,Hu2003}:
\begin{equation}
\label{eq.1}
\left\{\begin{array}{rlc}
\nabla \times \nabla \times u-k^2u+\nabla p=J&\text{in} &\Omega , \\
\nabla \cdot u=0&\text{in} &\Omega, \\
u\times n=0&\text{on}  &\partial \Omega , \\
p=0&\text{on} &\partial \Omega \\
\end{array}
\right.
\end{equation}
where $u$ is a vector field, $p$ is the scalar multiplier,
and $J$  is the given external source.
$\Omega$ is a simply connected domain in $\mathbb{R}^3$ with a connected boundary
$\partial \Omega$, with $n$ being its outward unit normal. The wave number~$ k$ is given by
$k^2=\omega^2 \varepsilon \mu,$ where $\omega $, $\varepsilon$ and $\mu$  are positive frequency, permittivity and
permeability of the medium, respectively. We assume that $k^2$ is not an interior Maxwell eigenvalue, and know
the cases with appropriately small and large frequencies are physically relevant in magnetostatics, wave propagation
and other applications \cite{Greif2007}. We refer to \cite[chapter 11]{Boffi2013} for a survey on this topic.
The introduction of the Lagrange multiplier $p$ in \eqref{eq.1} may not be absolutely necessary for the general case $k\ne 0$,
for which the divergence constraint does not need to be enforced explicitly; namely it is possible to
solve directly for $u$ using the first equation in \eqref{eq.1} with $p=0$ mathematically \cite{Hiptmair2002},
although it is still challenging to design an efficient numerical solver for this indefinite system.
The saddle-point formulation
\eqref{eq.1} with the Lagrange multiplier $p$
is stable and well-posed \cite{Demkowicz1998}, especially it ensures the stability and Gauss's law when $k$ is small
and may better handle the singularity of
the solution  at the boundary of the domain \cite{Demkowicz1998 , Perugia2002a,Chen:2000:FEM:588310.588552}.
More importantly, the mixed form \eqref{eq.1} provides some extra flexibility on the computational aspect \cite{Li2012}
and leads to better numerical stability
and more efficient numerical solvers than the single system \eqref{eq.1} without Lagrange multiplier (i.e., $p=0$),
as it was shown in \cite{Greif2015} \cite{Greif2007}. And this is also the main motivation and focus and of the current work.

After discretizing \eqref{eq.1} by using the N\'{e}d\'{e}lec elements of the first kind
\cite{Monk1992,Nedelec1980} for the approximation of the vector field $u$ and the standard nodal elements
for the multiplier $p$,  we derive the saddle-point system \eqref{eq.7} of our interest.
We  assume that the coefficient  matrix $\mathcal{K}$ in \eqref{eq.7} and its leading block $A-k^2M$ are both nonsingular, which is true
when the mesh size is sufficiently small \cite{Greif2007}.
%
%

Some very efficient preconditioners were proposed and analyzed recently in \cite{Greif2015} for
the special case of the saddle-point system \eqref{eq.7}, i.e.,
the wave number $k=0$, and
\eqref{eq.7} reduces to
\begin{equation}
\label{eq.2}
 \mathcal{A}\left(\begin{array}{c}u\\ p\end{array}\right)\equiv\begin{pmatrix}
A & B^T \\
B & 0
\end{pmatrix} \left(\begin{array}{c}u\\ p\end{array}\right)=\left(\begin{array}{c}f\\ g\end{array}\right)\,.
\end{equation}
%
The following preconditioner $\mathcal{P}_0$ was proposed in \cite{Greif2015}
for solving the saddle-point system~\eqref{eq.2}:
 \begin{equation}
 \label{eq.5}
\mathcal{P}^{-1}_0= \begin{pmatrix}
(A+M)^{-1}(I-B^TL^{-1}C^T) & CL^{-1} \\
L^{-1}C^T & 0
\end{pmatrix},
\end{equation}
where the matrix $L\in \mathbb{R} ^{m \times m}$ is the discrete Laplacian, while
$C\in \mathbb{R}^{n\times m}$ is a sparse matrix,
whose columns span $ker(A)$ and can be formed easily using the gradients of the standard nodal bases \cite{Greif2015}.
It was proved that the preconditioned system $\mathcal{P}_0^{-1}\mathcal{A}$ is
simply diagonal, given by
\begin{equation*}
\mathcal{P}_0^{-1}\mathcal{A}=\begin{pmatrix}
(A+ M)^{-1}(A+B^TL^{-1}B)&  0 \\
0 & I
\end{pmatrix}\,.
\end{equation*}
 As both $A+ M$ and $A+B^TL^{-1}B$ are symmetric positive definite, it allows us to apply a CG-like method
 for the preconditioned system $\mathcal{P}_0^{-1}A$ in a non-stardand inner product, even both $\mathcal{A}$ and $\mathcal{P}_0$ are indefinite.
For the more general case $k\ne 0$, the block tridiagonal preconditioners
\begin{equation}\label{eq:tri}
{\mathcal{M}}_{\eta,\varepsilon}=\begin{bmatrix}
A+(\eta-k^2)M & (1-\eta \varepsilon )B^T\\
0 & \varepsilon L
\end{bmatrix}
\end{equation}
with double variable relaxation parameters  $\eta>k^2$ and $\varepsilon\neq 0$
were studied in \cite{Greif2007,Zeng2012,Wu2013,Cheng2009}.

In this work, we construct some new preconditioners for \eqref{eq.2} and   \eqref{eq.7} respectively. We
will first show that for the case with vanishing wave number ($k=0$),
instead of the aforementioned efficient solver by using the preconditioner  \eqref{eq.5} or
\eqref{eq:tri} (with $k=0$), we can directly make use of the solution of the singular curl-curl system to construct
a more direct and efficient solver.

As it was shown in \cite{Greif2015} that preconditioners ${\cal P}_0^{-1}$ in \eqref{eq.5} work very effectively
for the special and simple case with vanishing wave number ($k=0$).
We shall demonstrate that similar
preconditioners can be constructed also for the saddle-point linear system \eqref{eq.7}
with more general and difficult cases, i.e., $k\neq0$, including high frequency waves.
And we will see analytically the spectral distributions of
these new preconditioners are quite similar to the ones of the existing effective
preconditioners \eqref{eq:tri}. But the new preconditioner can be applied with
CG iteration under a non-standard inner product although both the coefficient matrix ${\cal K}$ and
the new preconditioner are indefinite for $k\ne 0$, and numerically they will perform mostly better and more stable than
the existing preconditioners \eqref{eq:tri}.


The rest of the paper is arranged as follows.
We shall develop in Section~\ref{inverse} an important formula for computing the inverse of $\mathcal{K}$,
based on which we propose in Subsection~\ref{sec.31}
some more direct and efficient solver for the saddle-point system \eqref{eq.2} with vanishing wave number.
Then we shall propose a new preconditioner and compare its performance with existing preconditioners
for the saddle-point system~\eqref{eq.7} with general wave numbers,
and study and compare the spectral properties of the preconditioned matrices in Subsection~\ref{preconditioners}.
Numerical experiments are presented in Section~\ref{numerical}.

\section{Computing the inverse of $\mathcal{K}$}
\label{inverse}
\label{sec.2}
We shall derive in this section some formulae for computing the inverse of the matrix $\mathcal{K}$ in \eqref{eq.7}.
To do so, we first recall some useful properties of the matrices
$A$, $B$, $M$, $L$ and  $C$, which are introduced in the Introduction.
\begin{prop}
\label{pr.1}
The matrices $A$, $B$, $M$, $L$ and  $C$ have the following properties \cite{Greif2007,Greif2015}:
\begin{enumerate}[(i)]
\item  $L=BC$\,, \q
$AC=0$\,, \q
$MC=B^T$. \label{item.1}
\item  $ \mathbb{R}^n=ker(A)\oplus ker(B)$. \label{item.2}
\item  \label{item.3}
There exists a constant $\bar\alpha>0$ independent of mesh size such that
$u^TAu\geq \bar\alpha\,u^TMu\  $  ~$\forall\, u\in ker(B)$.
\item  \label{item.4}$u_A^TMu_B=0\ $ ~$\forall\,  u_A\in ker(A), \,u_B\in ker(B)$.
\item  \label{item.5} $u_A^TB^TL^{-1}Bu_A=u_A^TMu_A$  ~$\forall\, u_A\in ker(A)$.
\item \label{item.6} The inverse of $\mathcal{A}$ can be represented by
\begin{equation}\label{eq:inverseA}
\mathcal{A}^{-1}= \begin{pmatrix}
V & CL^{-1} \\
L^{-1}C^T & 0
\end{pmatrix},
\end{equation}
where  the diagonal block $V$ is given by
\begin{equation}
\label{eq.4}
V=(A+B^TL^{-1}B)^{-1}(I-B^TL^{-1}C^T)=(A+B^TL^{-1}B)^{-1}-CL^{-1}C^T\,.
\end{equation}
\end{enumerate}
 \end{prop}
 We give another proof for \eqref{item.6} of Proposition~\ref{pr.1} in the Appendix,  where  the  much more
 general cases are allowed, i.e. the cases when $\mathcal{A}$ is non-symmetric and its (2,2) block is nonzero.
\begin{lem}
For the matrix $V$ in \eqref{eq.4}, we have the following results:
\begin{align}
&VA=I-CL^{-1}B\,,  \q VB^T=0\,. \label{eq.8}\\
&AV=I-B^TL^{-1}C^T\,, \q BV=0\,.\label{eq.10}
\end{align}
\end{lem}

{\bf Proof.}
Using the fact that $L=BC$, it is easy to see from \eqref{eq.4} that
\begin{align*}
VB^T&=(A+B^TL^{-1}B)^{-1}(I-B^TL^{-1}C^T)B^T
=(A+B^TL^{-1}B)^{-1}(B^T-B^T)
=0\,.
\end{align*}
The first relation in \eqref{eq.8} follows readily from the fact that
 the (1,1) block of the following matrix
$$\mathcal{A}^{-1}\mathcal{A}=\begin{pmatrix}
V &CL^{-1}\\
L^{-1}C^T & 0
\end{pmatrix} \begin{pmatrix}
A & B^T \\
B & 0
\end{pmatrix}$$
 is an $n\times n$ identity matrix.
 Noting that  $\mathcal{A}$ is symmetric, we know $\mathcal{A}^{-1}$and  $V$ are both symmetric.
 Then the two identities in \eqref{eq.10} follow immediately by taking the transpose of both sides
 of each identity in \eqref{eq.8}.

Now we are ready to derive a formula for computing the inverse of
the matrix $\mathcal{K}$ in \eqref{eq.7}.  Recall that $\mathcal{K}$ and $A-k^2M$ are invertible.
We shall write the (1,1) block of the inverse $\mathcal{K}^{-1}$ as  $T$, then we have the following
representation of the inverse of the saddle-point matrix $\mathcal{K}$.
%
\begin{thm}\label{thm:K}
The inverse of $\mathcal{K}$ is given by
\begin{equation}
\label{eq.12}
\mathcal{K}^{-1}= \begin{pmatrix}
T & CL^{-1} \\
L^{-1}C^T & k^2 L^{-1}
\end{pmatrix},
\end{equation}
where  $T$ satisfies
\begin{equation}
\label{eq.13}
(A-k^2M)T=AV\,, \q  BT=0\,.
\end{equation}
\end{thm}

{\bf Proof.}
We write $\mathcal{K}^{-1}$ as a perturbation of $\mathcal{A}^{-1}$ in the form
\begin{equation}
\label{eq.15}
\mathcal{K}^{-1}= \mathcal{A}^{-1}+\begin{pmatrix}
X_1 & X_2 \\
X_3 & X_4
\end{pmatrix},
\end{equation}
then using the fact that $\mathcal{K}\mathcal{K}^{-1}=I$, namely
\begin{equation*}
\label{eq.16}
\left[\mathcal{A}+\begin{pmatrix}
-k^2M & 0 \\
0 & 0
\end{pmatrix}\right]
\cdot \left[\mathcal{A}^{-1}+\begin{pmatrix}
X_1 & X_2 \\
X_3 & X_4
\end{pmatrix}\right]=I\,,
\end{equation*}
we obtain by a direct computing that
\begin{equation}
\label{eq.17}
-k^2M(V+X_1)+AX_1+B^TX_3=0,
\end{equation}
\begin{equation}
\label{eq.18}
-k^2(B^TL^{-1}+MX_2)+AX_2+B^TX_4=0,
\end{equation}
\begin{equation}
\label{eq.19}
BX_1=0\,, \q BX_2=0.
\end{equation}
From \eqref{eq.15} we know that $V+X_1$ is the (1,1) block of $\mathcal{K}^{-1},$
so it follows from \eqref{eq.10} and \eqref{eq.19} that
\begin{equation*}
\label{eq.21}
BT=B(V+X_1)=0.
\end{equation*}
Multiplying  \eqref{eq.17} by $C^T$ we derive
$$-k^2B(V+X_1)+LX_3=0,$$
which gives
 \begin{equation}
 \label{eq.22}
 X_3=k^2L^{-1}B(V+X_1)=0.
 \end{equation}
Similarly,  multiplying  \eqref{eq.18}  by $C^T$ we obtain
$$-k^2(I+BX_2)+LX_4=0.$$  By combining this equality with  the second relation in \eqref{eq.19}, we
come to
 \begin{equation}
 \label{eq.23}
X_4=k^2L^{-1}.
 \end{equation}
Then we may substitute  \eqref{eq.23}  into \eqref{eq.18} to get
\begin{equation}
\label{eq.24}
(A-k^2M)X_2=0\,,
\end{equation}
which proves
$X_2=0$.

Noting that we have proved $X_3=0$,  then \eqref{eq.17} reduces to
$-k^2M(V+X_1)+AX_1=0,$ or   $(A-k^2M)(V+X_1)=AV$,
which completes the desired proof.

The following result is important to help us understand the leading block $T$ of the inverse of ${\cal K}$ in \eqref{eq.12}.
\begin{thm}
\label{th.inv}
The matrix $A+\eta B^TL^{-1}B-k^2M$ is non-singular for any $\eta\neq k^2$, and
its null space is exactly the same as that of $A$ for $\eta=k^2$.
\end{thm}
{\bf Proof.}
By means of
\eqref{item.2} of Proportion~\ref{pr.1}, we can write for any $u\in \mathbb{R}^n$ that
$u=u_A+u_B$ with $u_A\in ker(A)$ and $u_B\in ker(B)$. If $(A+\eta B^TL^{-1}B-k^2M)u=0,$ then
$(A-k^2M)u_B+\eta B^TL^{-1}Bu_A-k^2Mu_A=0$. As the columns of $C$ span the null space of $A$, there exists  $p\in\mathbb{R}^m$
such that
 $u_A=Cp$. So we see
 $(A-k^2M)u_B+(\eta  -k^2)B^Tp=0$. Multiplying its both sides by $C^T,$ we drive
 $p=0$, hence $(A-k^2M)u_B=0$, yielding that $u_B=0$. Hence we have proved $u=0,$ and also the non-singularity of
 the desired matrix.

Next, we consider the case with $\eta=k^2$. We show the two matrices $A+k^2 B^TL^{-1}B-k^2M$ and $A$ have the same null space.
First, we assume $u\in ker(A)$ and write
$u=u_A+u_B$ with $u_A\in ker(A)$ and $u_B\in ker(B)$,
then we see readily that $u=u_A=Cp$, hence $(A+\eta B^TL^{-1}B-k^2M)u=(\eta  -k^2)B^Tp=0$.

Now we assume $u$ is in the null space of $A+k^2 B^TL^{-1}B-k^2M$. We still write $u=u_A+u_B$,
and follow the earlier proof of the non-singularity of the matrix, but with $\eta=k^2$ now.
Then we shall deduce $(A-k^2M)u_B=0$, which implies $u_B=0$, hence we know $Au=0$.

The following result comes directly from \eqref{eq.13} and Theorem~\ref{th.inv}. And it introduces a very crucial
parameter $\eta$ to the expression of the leading block $T$ of the inverse of ${\cal K}$ in \eqref{eq.12},
and it can take an arbitrary value except for $\eta\neq k^2$.
\begin{coro}
\label{co.1}
For any $\eta\neq k^2$, it holds that
\begin{eqnarray}
T
&=&(A+\eta B^TL^{-1}B-k^2M)^{-1}(I-B^TL^{-1}C^T)\nb\\
&=&(A+\eta B^TL^{-1}B-k^2M)^{-1}-\frac{1}{\eta-k^2} CL^{-1}C^T.
\label{eq.26}
\end{eqnarray}
\end{coro}

Corollary~\ref{co.1} can be simplified for the special case $k=0$.
\begin{coro}
\label{co.3}
For any $\eta\ne 0$,  we have
\begin{equation}
V=(A+\eta B^TL^{-1}B)^{-1}(I-B^TL^{-1}C^T)=(A+\eta B^TL^{-1}B)^{-1}-\frac{1}{\eta} CL^{-1}C^T.
\label{eq:V}
\end{equation}
\end{coro}

It is very interesting to see from above that the two matrices $V $ and $T$ are independent of the
parameter~$\eta$, although their explicit representations in \eqref{eq.26} and \eqref{eq:V} look closely
depending on $\eta$.
In conclusion, we can easily see from  Theorem~\ref{thm:K} and Corollary~\ref{co.1} the following
formula with $\eta\neq k^2$ for computing the inverse of matrix ${\cal K}$ in \eqref{eq.7}, which forms the basis in our construction
of some new preconditioners that are discussed in the next section:
\begin{equation}
\label{eq.27}
\mathcal{K}^{-1}= \begin{pmatrix}
(A+\eta B^TL^{-1}B-k^2M)^{-1}(I-B^TL^{-1}C^T) & CL^{-1} \\
L^{-1}C^T & k^2 L^{-1}
\end{pmatrix}.
\end{equation}

\section{New preconditioners and their spectral properties}
\label{sec.3}
\subsection{Vanishing wave number: $k=0$}
\label{sec.31}
In this subsection we will propose a new simple solution to the saddle-point system \eqref{eq.2}.
To do this, we can easily see from Proporsition\,\ref{pr.1} (vi) that $p=L^{-1}C^T f$, then we can rewrite
the saddle-point system \eqref{eq.2} as
\begin{eqnarray}
Au &=& (I-B^TL^{-1}C^T)f\,, \label{eq.k1} \\
Bu &=&g\,.   \label{eq.k2}
\end{eqnarray}
%
%
As $A$ is the discrete counterpart of the curl\,curl operator, so it is singular and the equation \eqref{eq.k1} has multiple solutions.
Now we consider an arbitrary solution $u_0$ to \eqref{eq.k1}. Noting that the columns of $C$ span the kernel of $A$ (see section 1),
we can decompose the solution $u$ of
\eqref{eq.2} into $u=u_0+Ct$ for some vector $t\in \mathbb{R}^m$.
By means of \eqref{eq.k2} and the fact that $BC=L$ from Proporsition\,\ref{pr.1}, we derive
from $Bu=B(u_0+Ct)=g$ that $t=L^{-1}(g-Bu_0)$, which yields
$$u=(I-CL^{-1}B)u_0+CL^{-1}g,$$
or
\begin{equation}
\label{eq.k3}
 \left(\begin{array}{c}u\\ p\end{array}\right)
=\begin{pmatrix}
(I-CL^{-1}B)  A^{+}(I-B^TL^{-1}C^T)  & CL^{-1} \\
L^{-1}C^T & 0
\end{pmatrix} \left(\begin{array}{c}f\\ g\end{array}\right),
\end{equation}
which gives another representation of the inverse $\mathcal{A}^{-1}$ and its diagonal part $V$ (see \eqref{eq:inverseA}).
Here $A^{+}$ is any operator which maps any vector $b\in \mathcal{R}(A)$ to a particular solution
of $Au=b$.

We remark that if the source $J$ is divergence free in \eqref{eq.1},  then the Lagrange multiplier $p=0$, so is the discrete
$p$ in the saddle-point system \eqref{eq.2}. Then we get from $p=L^{-1}C^T f$ that $C^Tf=0$. In this case, the equation
\eqref{eq.k1} reduces to $Au=f$. We have similar simplification in \eqref{eq.k3}.


To find an arbitrary solution $u_0$ to \eqref{eq.k1},
we may apply the CG iteration with the Hiptmair-Xu preconditioner \cite{Hiptmair2007}, which works
vey efficiently for the discrete system \eqref{eq.k1} arising from the discretization of the curl curl system \cite{rep}.
One may also develop a preconditioner for the whole system based on the important equation
\eqref{eq.k3}.

Comparing with the preconditioner ${\cal P}_0^{-1}$ in \eqref{eq.5} or the preconditioner in \eqref{eq:tri} ($k=0$),
the above new solver generated by the relation \eqref{eq.k3} should be much more efficient computationally
as the new solver may be viewed like a direct solver.

\subsection{General wave numbers: $k\ne 0$}
\label{preconditioners}
The formula \eqref{eq.27} suggests us some natural
preconditioners for the saddle-point matrix ${\cal K}$ in \eqref{eq.7}.
Noting that the matrix $B^TL^{-1}B$ is a dense matrix,
the action of the (1,1) block of \eqref{eq.27} is very expensive to compute.
To overcome the difficulty, we approximate the dense matrix $A+\eta B^TL^{-1}B-k^2M$
by the sparse matrix $A+\eta M-k^2M$ to get
the following simplified preconditioner for the matrix ${\cal K}$:
\begin{equation}
\label{eq.28}
\mathcal{P}^{-1}\equiv\begin{pmatrix}
(A+\eta M-k^2M)^{-1}(I-B^TL^{-1}C^T) & CL^{-1} \\
L^{-1}C^T & k^2 L^{-1}
\end{pmatrix}.
\end{equation}
For the simple case with vanishing wave number ($k=0$) and
$\eta=1$, the preconditioner \eqref{eq.28} reduces to the existing one $\mathcal{P}^{-1}_0$ in \eqref{eq.5}.
To ensure the nonsingularity of the matrix $A+\eta M-k^2M$ involved in \eqref{eq.28}, we can simply
set the parameter $\eta>k^2$ so that it becomes symmetric positive definite.
Moreover, this choice also ensures the nonsingularity of the matrix on the right-hand side of \eqref{eq.28},
as discussed below.
\begin{thm}
For any $\eta>k^2$, the matrix on the right-hand side of \eqref{eq.28} is non-singular.
\end{thm}

{\bf Proof.}
It is direct to check that the preconditioned matrix $\mathcal{P}^{-1}\mathcal{K}$ is given by
\begin{equation*}
\label{eq.29}
\mathcal{P}^{-1}\mathcal{K}=\begin{pmatrix}
(A+\eta M-k^2M)^{-1}(A+k^2B^TL^{-1}B -k^2M)+CL^{-1}B & 0 \\
0 & I
\end{pmatrix}\,.
\end{equation*}
Using Proposition~\ref{pr.1}\,\eqref{item.1}, we can further write the
(1, 1) block of the above matrix as follows:
\begin{equation*}
\begin{split}
&(A+\eta M-k^2M)^{-1}(A+k^2B^TL^{-1}B -k^2M)+CL^{-1}B\\
=&(A+\eta M-k^2M)^{-1}(A+k^2B^TL^{-1}B -k^2M)\\
&+(A+\eta M-k^2M)^{-1}(\eta MCL^{-1}B-k^2 MCL^{-1}B)\\
=&(A+\eta M-k^2M)^{-1}(A+\eta B^TL^{-1}B-k^2M)\,,
\end{split}
\end{equation*}
so the preconditioned matrix $\mathcal{P}^{-1}\mathcal{K}$ reads also as
\begin{equation}
\label{eq.30}
\mathcal{P}^{-1}\mathcal{K}=\begin{pmatrix}
(A+\eta M-k^2M)^{-1}(A+\eta B^TL^{-1}B-k^2M) & 0 \\
0 & I
\end{pmatrix}.
\end{equation}
We know that the leading block of $\mathcal{P}^{-1}\mathcal{K}$ in \eqref{eq.30} is non-singular
by Theorem~\ref{th.inv}, hence the desired conclusion follows.

Note that $A+\eta M-k^2M$ and its inverse are always symmetric positive definite for $\eta>k^2$.
Actually, the original matrix $A+\eta B^TL^{-1}B -k^2M$ can be also symmetric positive definite
as shown below.
\begin{thm}
\label{tm.bound}
If any $\eta>k^2$ and $k^2< \bar\alpha$, the matrix $A+\eta B^TL^{-1}B -k^2M$ is symmetric positive definite.
\end{thm}

{\bf Proof.}
For any $u\in \mathbb{R}^n$, we can write $u=u_A+u_B$ with $u_A\in ker(A)$ and $u_B\in ker(B)$.
By  Proposition~\ref{pr.1}  we know $u_A^TMu_B=0$ and $u_A^TB^TL^{-1}Bu_A=u_A^TMu_A$. Therefore,
\begin{equation}
\label{eq.long}
 \begin{split}
&u^T(A+\eta B^TL^{-1}B -k^2M)u\\
=&u_A^T(A+\eta B^TL^{-1}B -k^2M)u_A+u_B^T(A+\eta B^TL^{-1}B -k^2M)u_B\\
=&u_A^T(\eta B^TL^{-1}B -k^2M)u_A+u_B^T(A -k^2M)u_B\\
=&u_B^T(A-k^2M)u_B+(\eta-k^2) u_A^TMu_A\,.
 \end{split}
\end{equation}
But we have $u_B^TAu_B\geq\bar\alpha u_B^TMu_B$ by \eqref{item.3}  in Proposition~\ref{pr.1}, and this
implies
$$u^T(A+\eta B^TL^{-1}B -k^2M)u\geq (\bar\alpha-k^2) u_B^TMu_B+(\eta-k^2) u_A^TMu_A>0\,,$$
hence proves the desired result.


For $\eta>k^2$ and $k^2< \bar\alpha$, the preconditioned  matrix $
\mathcal{P}^{-1}\mathcal{K}$ is self-adjoint and positive definite with respect to the inner product
\begin{equation}
\label{inpro}
\left<x, y \right>=x^T\begin{pmatrix}
(A+\eta M-k^2M) & 0 \\
0 & I
\end{pmatrix}y. \end{equation}
So we can apply the CG iteration \cite{Ashby1990} in this special inner product for solving
the preconditioned system associated with ${\cal P}^{-1}{\cal K}$.
But Theorem~\ref{tm.bound} does not ensure the  positive definiteness
of the preconditioned system for $k^2\geq\bar\alpha$,  so the CG iteration may fail theoretically. In this case,
we may still apply the preconditioned MINRES with the above non-stardand inner product.

We know the convergence rates of the CG and MINRES can be reflected
often by the spectrum of the preconditioned system.
For this purpose, we shall now study the spectral properties of the preconditioned system ${\cal P}^{-1}{\cal K}$.
First, we present an interesting observation that
the parameter~$\eta$ does not affact the symmetric positive definiteness of the matrix $A+\eta  B^TL^{-1}B -k^2M$.
\begin{thm}
\label{th.negative}
For any two numbers $\eta_1,\eta_2>k^2,$
$A+\eta_1 B^TL^{-1}B -k^2M$ is  symmetric positive definite if and only if  $A+\eta_2 B^TL^{-1}B -k^2M$ is   symmetric positive definite.
\end{thm}
{\bf Proof.}
For any $\eta_1>k^2,$ suppose that $A+\eta_1 B^TL^{-1}B -k^2M$ is not symmetric positive definite.
As this matrix is nonsingular by Theorem~\ref{th.inv},   hence it is not symmetric semi-positive definite.
Therefore, there exists $u\in \mathbb{R}^n$ satisfying $u^T(A+\eta_1 B^TL^{-1}B -k^2M)u<0$.
But we can write $u=u_{A}+u_{B}$ with $u_{A}\in ker(A)$ and $u_{B}\in ker(B)$.
Then we can see that $u_B\neq 0$ and $u_{B}^T(A-k^2M)u_{B}<0$ from \eqref{eq.long}.
Now for any $\eta_2>k^2,$ we can easily check
 $u_B^T(A+\eta_2 B^TL^{-1}B -k^2M)u_B=u_{B}^T(A-k^2M)u_{B}<0$, hence
$A+\eta_2 B^TL^{-1}B -k^2M$ is not symmetric positive definite.

Next we present several results about the eigenvalues of the preconditioned matrix
$\mathcal{P}^{-1}\mathcal{K}$.
\begin{lem}
\label{lem.1}
For  any $\eta>k^2$, $\lambda =1$ is  an eigenvalue of
$(A+\eta M-k^2M)^{-1}(A+\eta B^TL^{-1}B -k^2M)$
with its algebraic multiplicity being $m$. The rest  of the eigenvalues are bounded by
\begin{equation}
\label{eq.31}
\frac{\bar\alpha-k^2}{\bar\alpha+\eta-k^2}<\lambda<1 .
\end{equation}
\end{lem}

{\bf Proof.} The result was proved in \cite[Theorem~5.1]{Greif2007} for $\eta=1$ and $k^2<1$.
But our desired results for an arbitrary positive $\eta$ can be done similarly.

The following result is a direct consequence of Lemma \ref{lem.1} by using the formula \eqref{eq.30}.
\begin{thm}
\label{th.6}
For  any $\eta>k^2$,   $\lambda =1$ is  an eigenvalue of
the preconditioned matrix $\mathcal{P}^{-1}\mathcal{K}$
with its algebraic multiplicity being $2m$. The rest  of the eigenvalues are bounded as in \eqref{eq.31}.
\end{thm}

Now we like to make some spectral comparisons between
the two preconditioned systems generated by our new preconditioner
$\mathcal{P}$  and the existing block triadiagonal one $\mathcal{M}_{\eta,\varepsilon}$ in \eqref{eq:tri}
for the saddle-point matrix $\mathcal{K}$.
%
We first recall the following results from \cite[Theorem 2.6]{Zeng2012}.
\begin{thm}
\label{th.7}
For any $\eta>k^2$,  both $\lambda_1=1 $ and $\lambda_2=-\frac{1}{\varepsilon(\eta-k^2) }$
are the eigenvalues of $\mathcal{M}_{\eta,\varepsilon}^{-1}\mathcal{K},$ each with
its algebraic multiplicity $m$.  The rest  of the eigenvalues are bounded as in \eqref{eq.31}.
\end{thm}
%

We see from Theorems~\ref{th.6} and \ref{th.7} that
the spectra of $\mathcal{P}^{-1}\mathcal{K}$ and  $\mathcal{M}_{\eta,\varepsilon}^{-1}\mathcal{K}$  are quite similar,
except that the latter has an extra eigenvalue $\lambda_2$, with its algebraic multiplicity being $m$.
This will be also confirmed numerically in the next section.


The block triadiagonal preconitioners $\mathcal{M}_{\eta,\varepsilon}$ reduce
to symmetric if we set $\varepsilon=\frac{1}{\eta}$:
\begin{equation}\label{eq:diagonal}
\mathcal{M}_{\eta, 1/\eta}=\begin{bmatrix}
A+(\eta-k^2)M & 0\\
0 & \frac{1}{\eta} L
\end{bmatrix}\,.
\end{equation}
This preconditioner was analysed and applied in \cite{Greif2007, Wu2013} along with the MINRES iteration.
We may observe from Theorems~\ref{th.6}  and \ref{th.7}
that the eigenvalues of our preconditioned matrix  $\mathcal{P}^{-1}\mathcal{K}$
are a little better clustered  than those of $\mathcal{M}_{\eta, 1/\eta}^{-1}\mathcal{K}$
as its eigenvalue $\lambda_2$ is smaller than $\frac{\bar\alpha-k^2}{\bar\alpha+\eta-k^2} $.
But our new preconditioner $\mathcal{P} $ can be applied with CG
for $k^2<\bar\alpha,$ and MINRES for $k^2\geq\bar\alpha$.
And more importantly, as we will see from our numerical experiments in next section,
we can also apply the new preconditioner $\mathcal{P} $ with CG even for
$k^2\geq\bar\alpha$ and the convergence is still rather stable,
while CG with preconditioner  $\mathcal{M}_{\eta,\varepsilon}$ in \eqref{eq:diagonal} breaks down most of the time.


On the other hand, if we choose
$\varepsilon\ne {1}/{\eta}$, the preconditioner  $\mathcal{M}_{\eta,\varepsilon}$ is non-symmetric,
and the methods like GMRES should be used, which are less economical than methods like CG or MINRES.
Note that for $\varepsilon=-\frac{1}{\eta-k^2},$ we have $\lambda_2=\lambda_1$, so $\lambda =1$ is
an eigenvalue of $\mathcal{M}_{\eta,\varepsilon}^{-1}\mathcal{K}$
with  its algebraic multiplicity being $2m,$ the same as for $\mathcal{P}^{-1}\mathcal{K}$.



Now we consider the inner iterations associated with the new preconditioner $\mathcal{P}$.
For any two vectors $ x\in\mathbb{R}^{n}$ and $y\in\mathbb{R}^{m}$,
we can write
\begin{equation*}
\begin{split}
\mathcal{P}^{-1} \left(\begin{array}{c}x\\ y\end{array}\right)=&
\left(\begin{array}{c}(A+\eta M-k^2M)^{-1}(x-B^TL^{-1}C^Tx)+CL^{-1}y \\
L^{-1}C^Tx+ k^2 L^{-1}y\end{array}\right)\, \\
=&\left(\begin{array}{c}(A+\eta M-k^2M)^{-1}x-\frac{1}{\eta-k^2}CL^{-1}C^Tx+CL^{-1}y \\
	L^{-1}C^Tx+ k^2 L^{-1}y\end{array}\right)\,.
\end{split}
\end{equation*}
So we have to solve two linear systems associated with  the discrete Laplacian $L$ and one with $A +(\eta-k^2) M$ at each evaluation of the action of $\mathcal{P}^{-1}$.
Many fast solvers are available for solving these two symmetric and positive definite systems
\cite{Li2012,Hiptmair2007}.
We know from Theorem~\ref{th.6} that a small $\eta-k^2$ may result in a better convergence of
the preconditioned Krylov subspace methods. But  if  $\eta-k^2$ is  too small,
the matrix $A +(\eta-k^2) M$ would become nearly singular.

We know that the parameter~$\bar\alpha $  depends only on the
shape regularity of the mesh and the approximation order of the finite elements used,
and is irrelevant to the  size of the mesh \cite[Theorem 4.7]{Hiptmair2002}.
Numerically we may expect a upper bound for $k^2$  that ensures
the positive definiteness of $A+\eta B^TL^{-1}B-k^2M$, and  this bound should be independent of the mesh size.
We shall test this numerically in the next section.

\subsection{Preconditioners for more general saddle-point systems}
We devote this subsection to discuss a by-product produced by our previous technique
for deriving the formula \eqref{eq.27} to compute the inverse ${\cal K}^{-1}$,
yielding an effective preconditioner for the following saddle-point matrix
with an arbitrary and non-vanishing (2, 2) block $D\in\mathbb{R}^{m\times m}$:
$$
\mathcal{A}_D\equiv\begin{pmatrix}
A & B^{T} \\
B & D
\end{pmatrix}\,.
$$

In fact, we can easily check using   Proposition~\ref{pr.1} that
\begin{equation}
\label{sp.ad}
\begin{pmatrix}
A & B^{T} \\
B & D
\end{pmatrix}\begin{pmatrix}
V-CL^{-1}DL^{-1}C^T & CL^{-1} \\
L^{-1}C^T & 0
\end{pmatrix}=I\,,
\end{equation}
so  $\mathcal{A}_D$ is non-singular, and its inverse is given by the second matrix
on the left-hand side of \eqref{sp.ad}.
Formula \eqref{sp.ad} is very
interesting: the inverse of $\mathcal{A}_D$ does not involve
the inverse of matrix $D$.
Note that \eqref{sp.ad} holds not only for those systems arising from Maxwell equations
but also for all non-singular saddle-point matrices with their $(2,2)$  block  being nonzero and
their leading block being maximally
rank deficient; see \cite{Greif2015} for more applications.

Motivated by \eqref{eq.28}, we may consider the following preconditioner:
\begin{equation}
\mathcal{P}_D^{-1}=\begin{pmatrix}
(A+\eta M)^{-1}(I-B^TL^{-1}C^T)-CL^{-1}DL^{-1}C^T & CL^{-1} \\
L^{-1}C^T & 0
\end{pmatrix}\,.
\end{equation}
Then we readily see a block diagonal preconditioned system:
\begin{equation}
{\mathcal{P}_D}^{-1}\mathcal{A}_D
=\begin{pmatrix}
(A+\eta M)^{-1}(A+\eta B^TL^{-1}B) & 0 \\
0 & I
\end{pmatrix}.
\end{equation}
This important relation indicates that we can always apply the CG iteration for solving the saddle-point system
associated with the matrix ${\cal A}_D$
in a special inner product as we did in section\,\ref{preconditioners}.

We can easily see that for the special case that $D=-\frac{1}{\eta}L$,  preconditioner $\mathcal{P}_D^{-1}$ reduces to
\begin{equation*}
\label{sp.pd}
\mathcal{P}_D^{-1}=\begin{pmatrix}
(A+\eta M)^{-1} & CL^{-1} \\
L^{-1}C^T & 0
\end{pmatrix}.
\end{equation*}

%

\section{Numerical experiments}
\label{numerical}
\begin{figure}[!hbt]\scriptsize
\caption{Meshes G1 through G4}
\label{fi.1}

\centering

\subfigure[\scriptsize{Mesh G1: n+m=187}]{\begin{minipage}[htb]{0.25\textwidth}
\label{Fi.g1}
\includegraphics[width=1\textwidth]{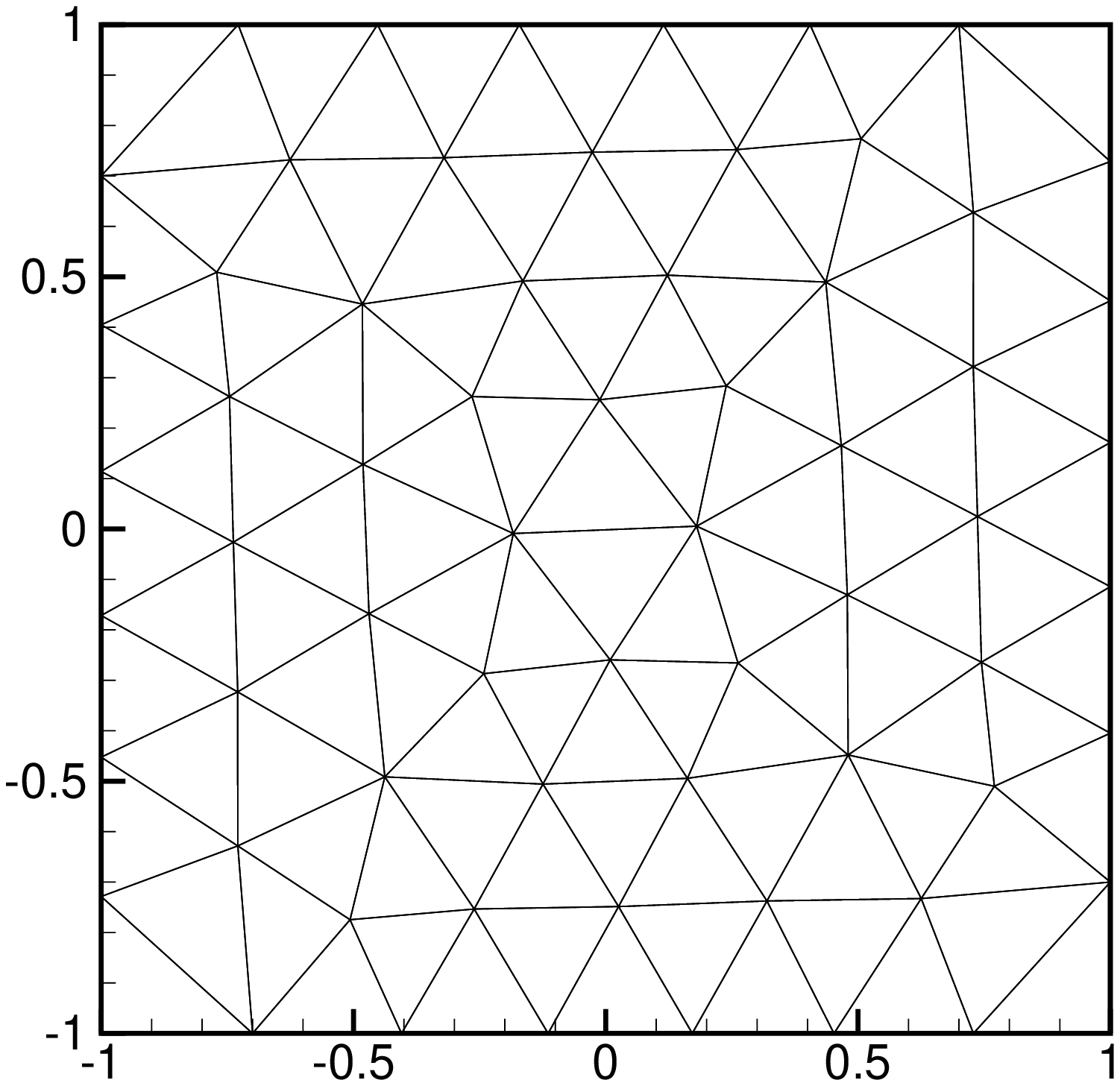}\end{minipage}}
\subfigure[\scriptsize{Mesh G2: n+m=437}]{\begin{minipage}[htb]{0.25\textwidth}
\label{Fi.g2}
\includegraphics[width=1\textwidth]{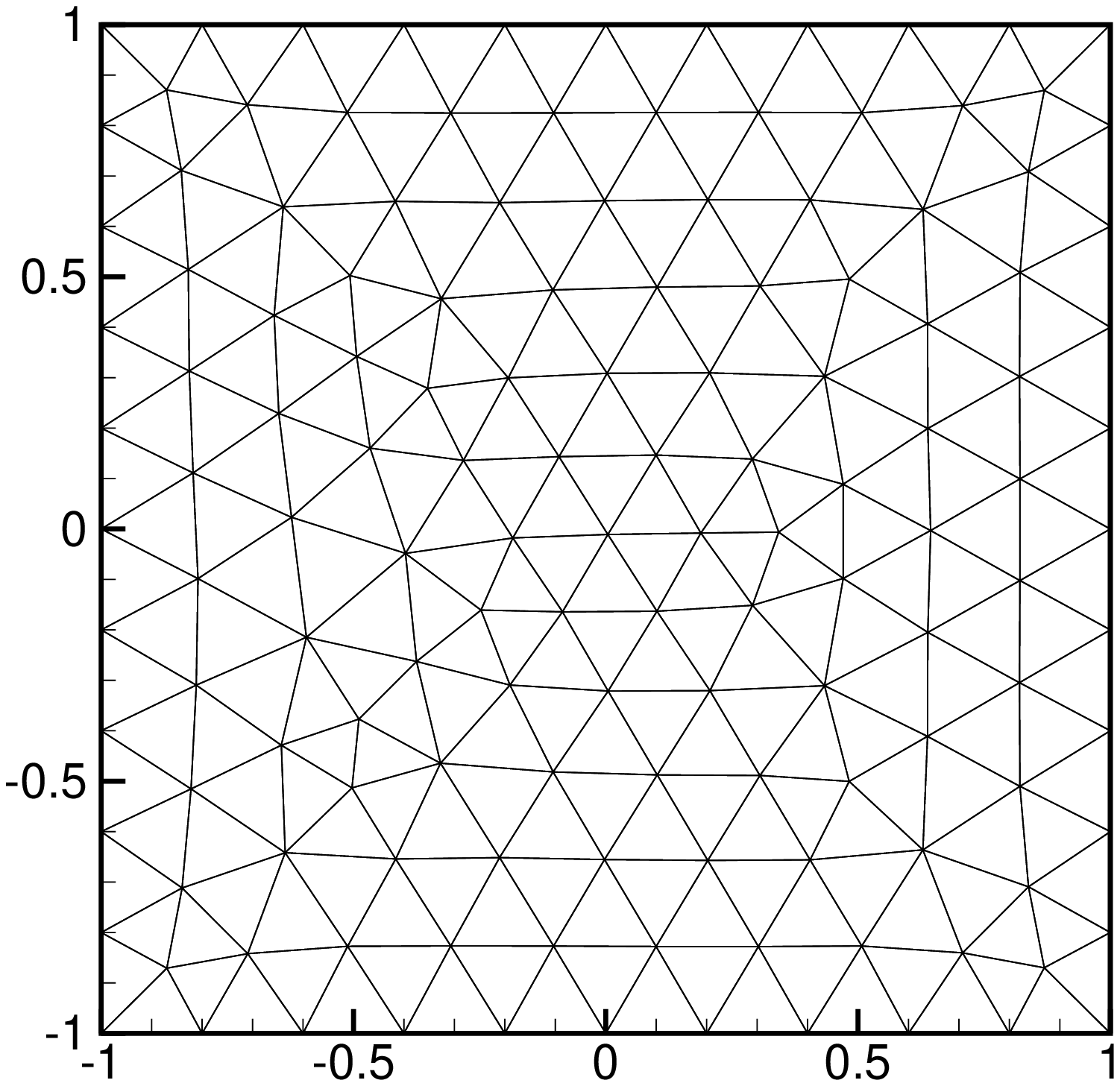}\end{minipage}}

\subfigure[\scriptsize{Mesh G3: n+m=1777}]{\begin{minipage}[htb]{0.25\textwidth}
\label{Fi.g3}
\includegraphics[width=1\textwidth]{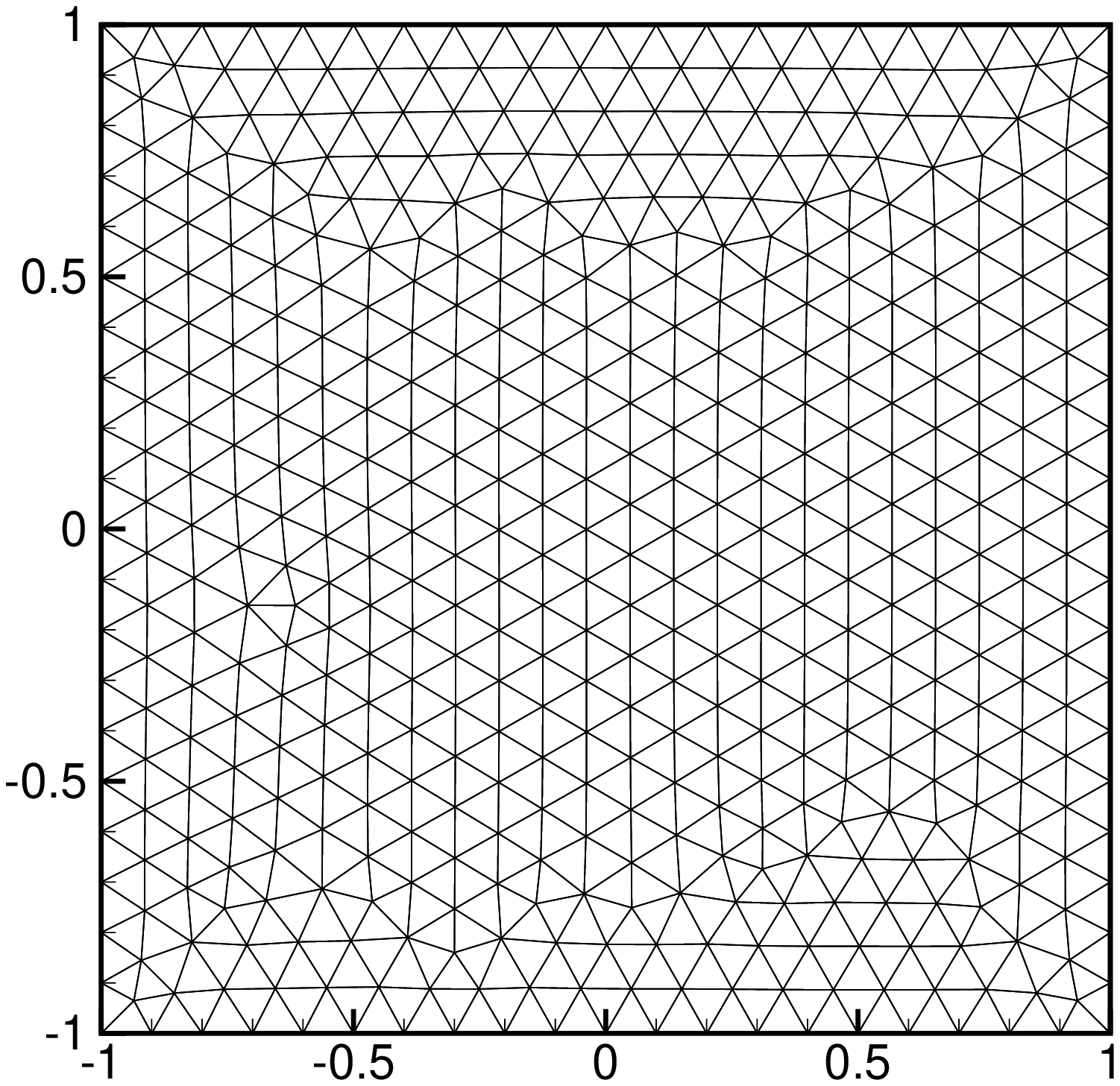}\end{minipage}}
\subfigure[\scriptsize{Mesh G4: n+m=7217}]{\begin{minipage}[htb]{0.25\textwidth}
\label{Fi.g4}
\includegraphics[width=1\textwidth]{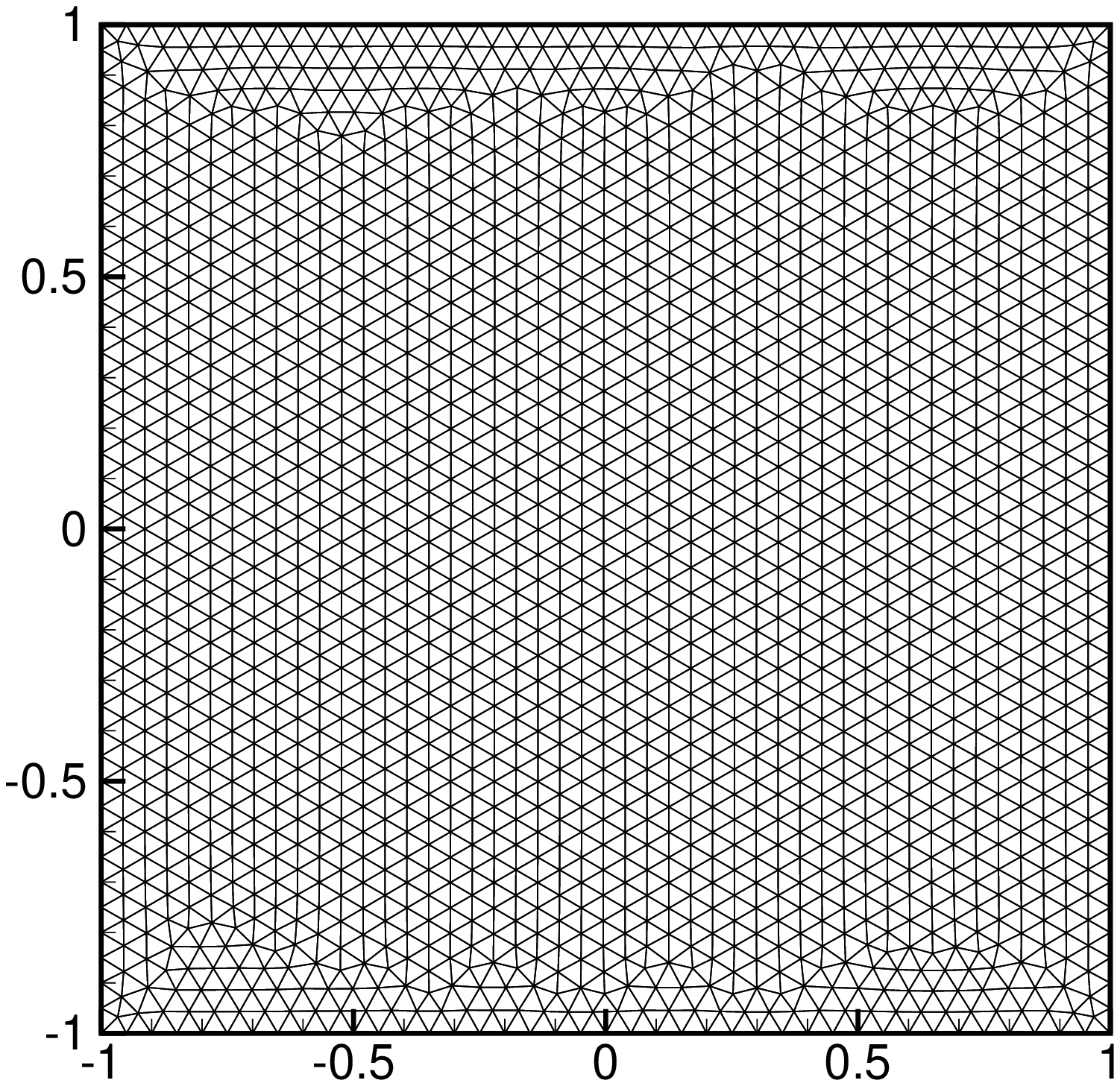}\end{minipage}}

\end{figure}

\begin{figure}[!hbt]\scriptsize
\vspace{.7cm}
\caption{Grids L1 through L4}
\label{fi.2}

\centering

\subfigure[\scriptsize{L1, with n+m=185}]{\begin{minipage}[htb]{0.25\textwidth}
\label{Fi.l1}
\includegraphics[width=1\textwidth]{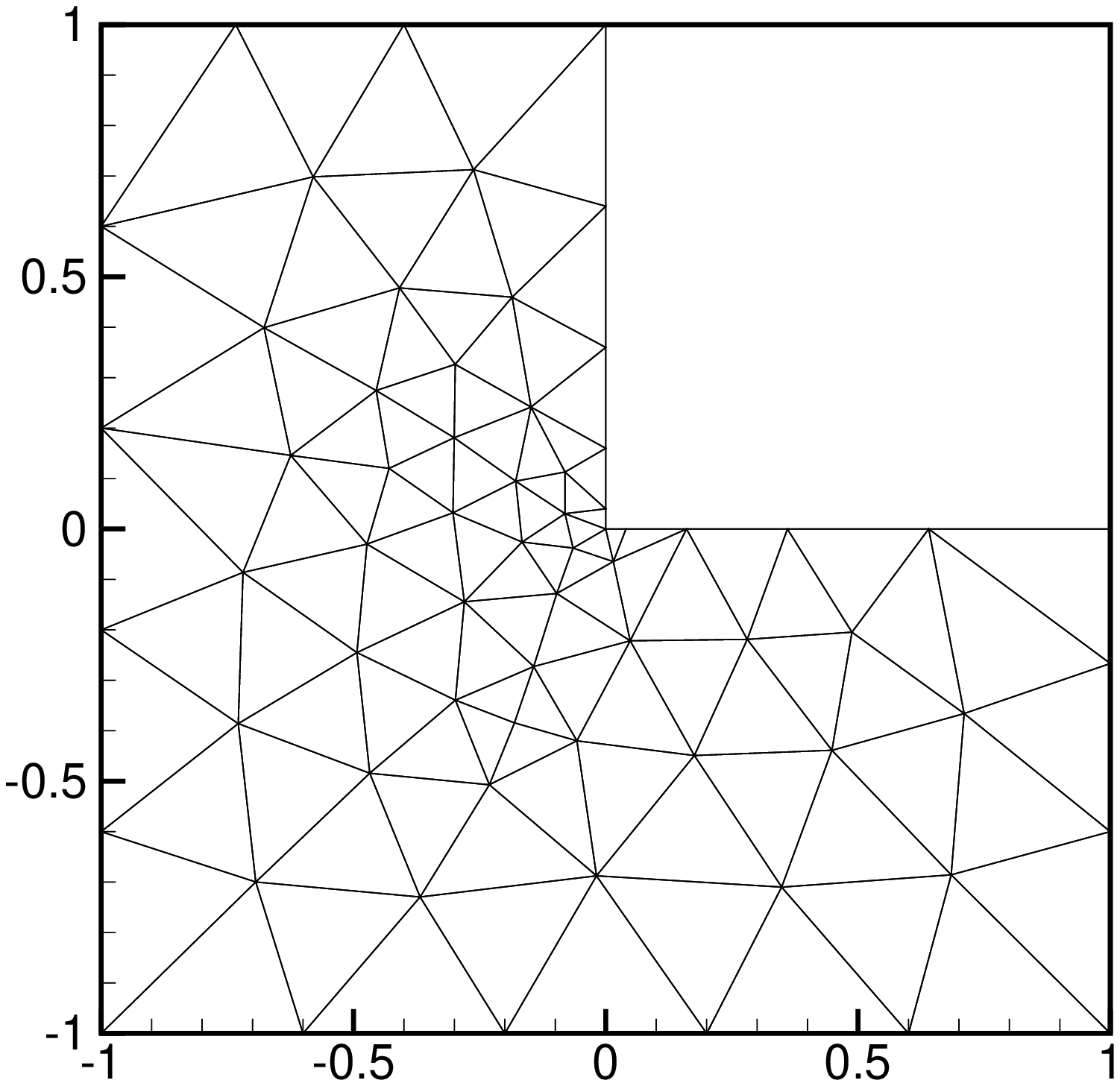}\end{minipage}}
\subfigure[\scriptsize{L2, with n+m=409}]{\begin{minipage}[htb]{0.25\textwidth}
\label{Fi.l2}
\includegraphics[width=1\textwidth]{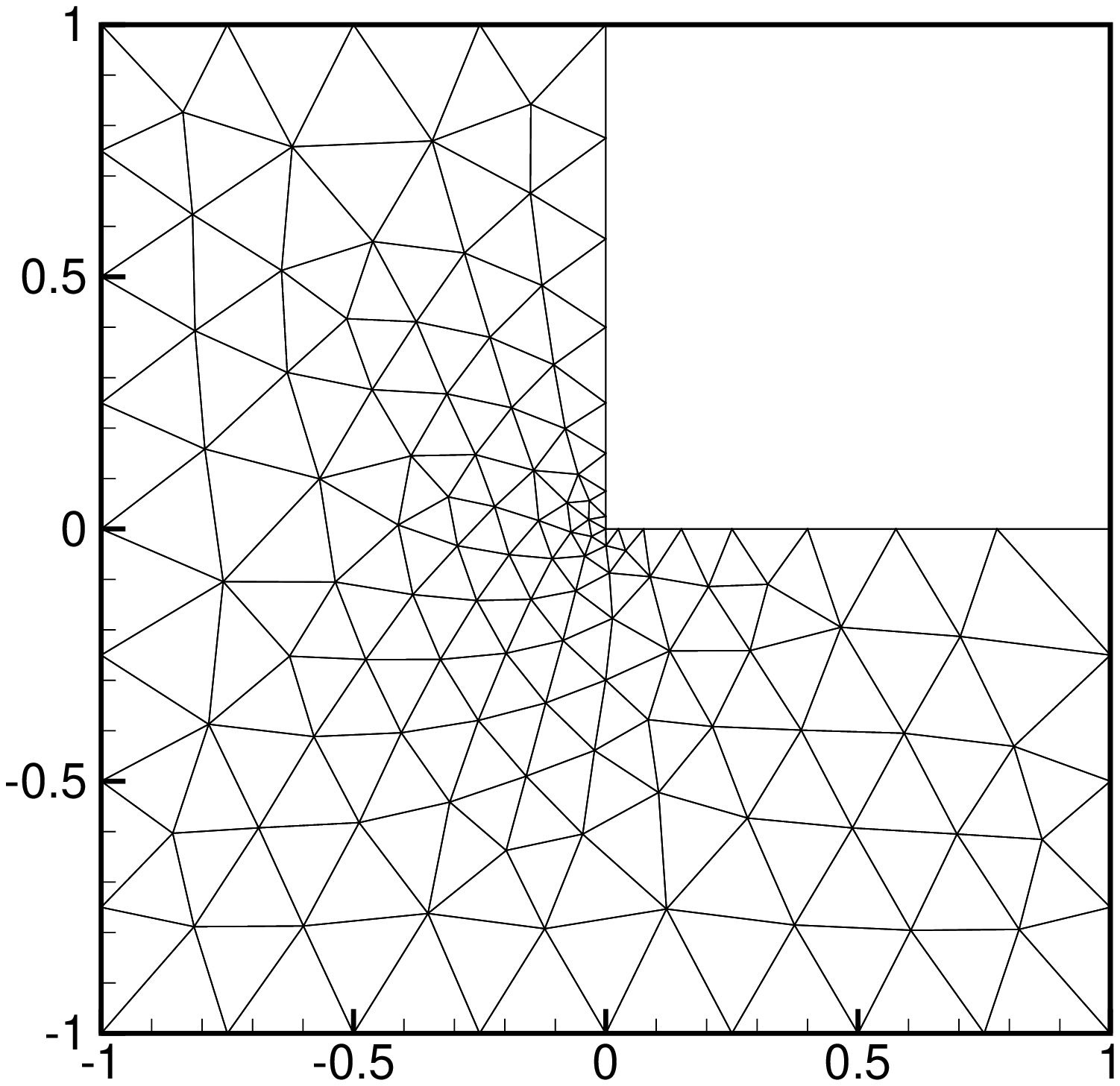}\end{minipage}}

\subfigure[\scriptsize{L3, with n+m=1177}]{\begin{minipage}[htb]{0.25\textwidth}
\label{Fi.l3}
\includegraphics[width=1\textwidth]{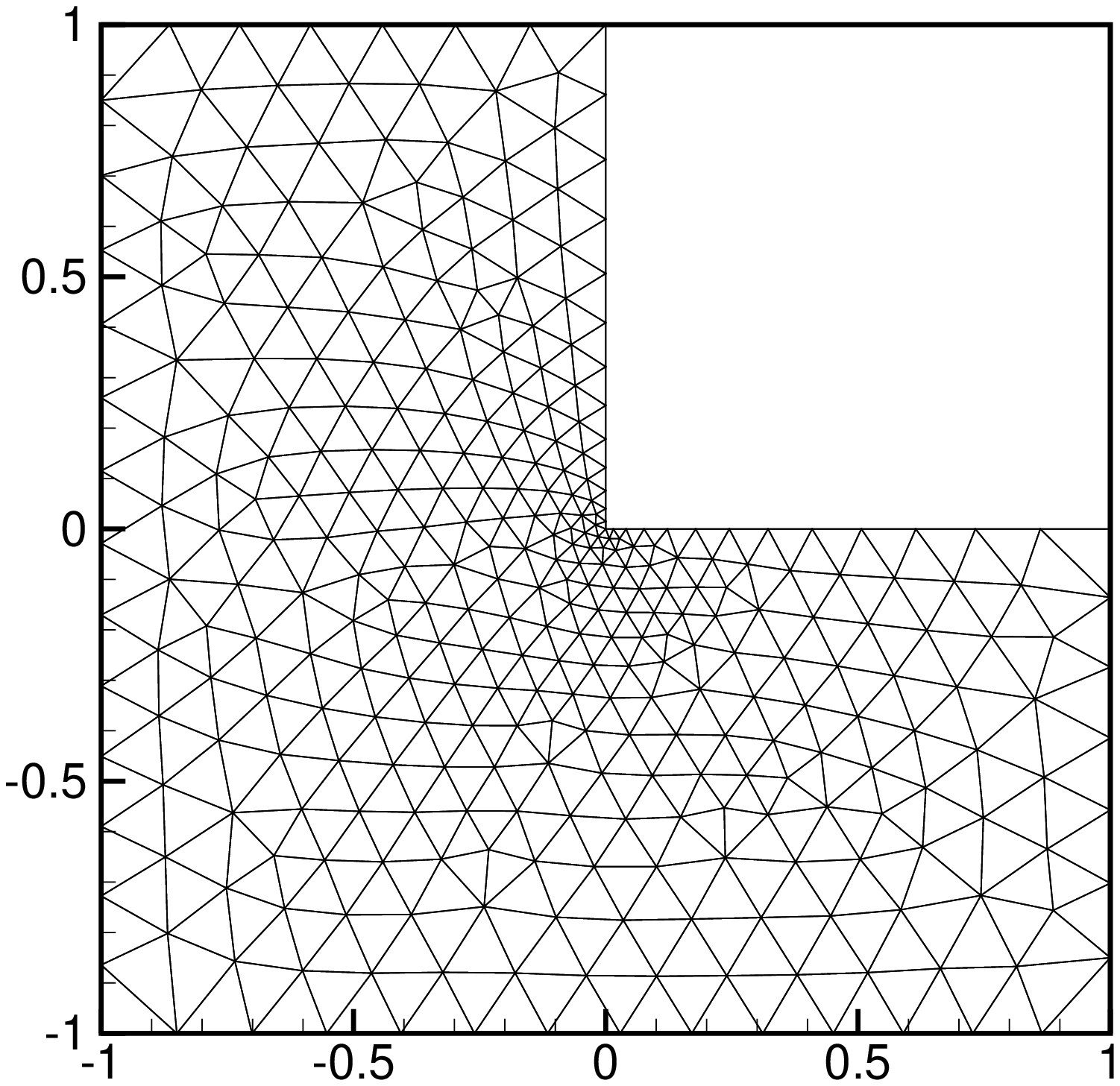}\end{minipage}}
\subfigure[\scriptsize{L4, with n+m=5325}]{\begin{minipage}[htb]{0.25\textwidth}
\label{Fi.l4}
\includegraphics[width=1\textwidth]{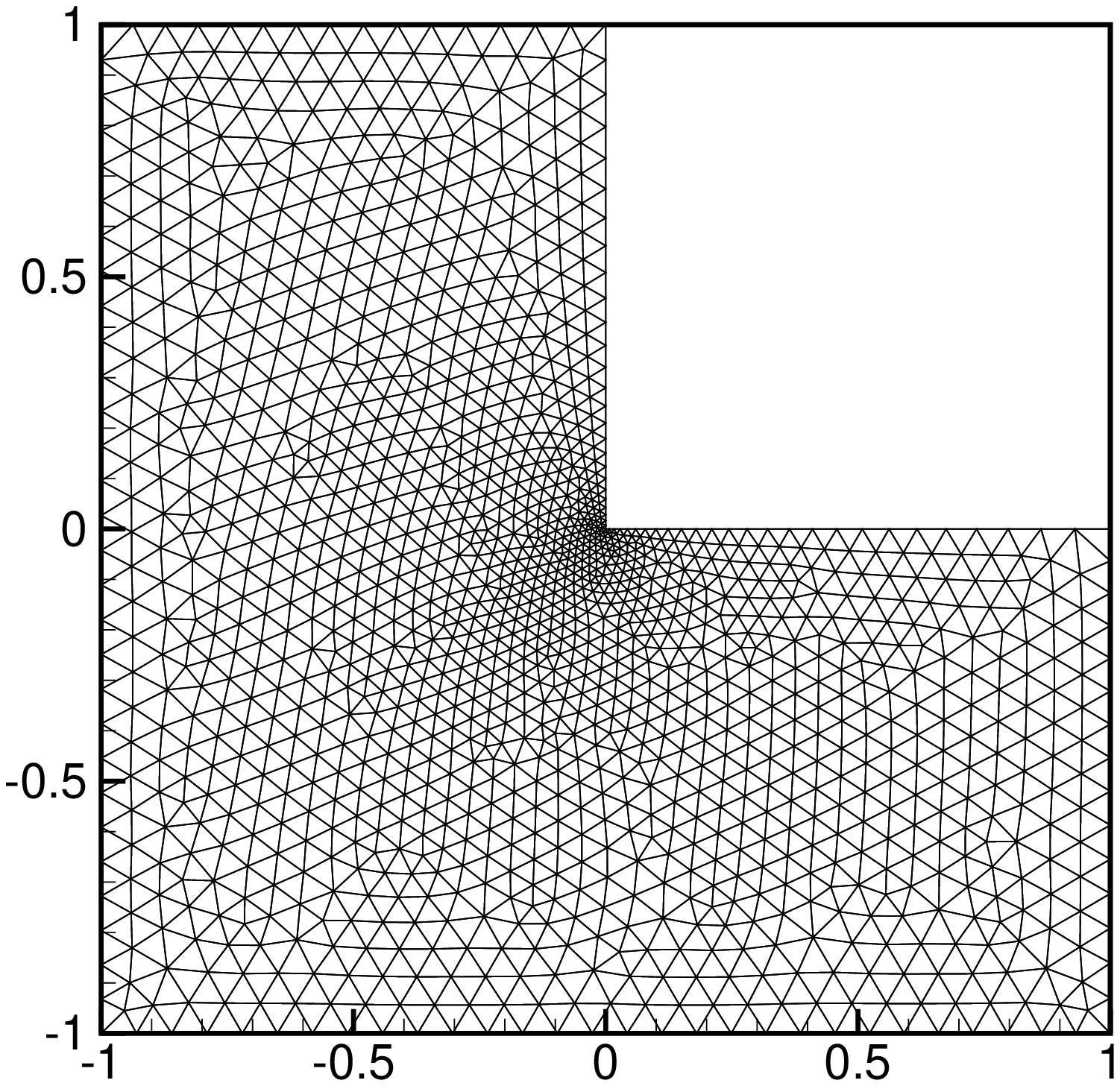}\end{minipage}}
\end{figure}

In this section, we present numerical experiments to demonstrate and compare
the spectral distributions of
the preconditioned systems of the saddle-point problem \eqref{eq.7} with the existing
preconditioner ${\mathcal{M}}_{\eta, 1/\eta}$ in  \eqref{eq:tri} and the new one
${\cal P}$ in \eqref{eq.28}, and the results of some Krylov subspace  iterations.
The  edge elements of lowest order  are used for the discretization of the system \eqref{eq.1}
in a square domain $(-1 \leq x \leq 1, -1 \leq y \leq 1 )$ or an L-shaped domain (see Figures~\ref{fi.1} and \ref{fi.2}),
which is partitioned using unstructured simplicial meshes generated
by EasyMesh \cite{D}.  For Meshes G1 through G5,  the desired side lengths of the triangles that contain one of the vertices of the domain 
are set to be the same.
For Meshes L1 through L5, the desired side lengths of the triangles that contain the  origin are one-tenth of the desired side lengths of
the triangles that contain other vertices of the domain.
Meshes G5 and L5  lead to linear systems of size $n+m=23769$ and $29277$ respectively.
We use \textsc{Matlab}  to implement all numerical iterative solvers,
and the linear solvers involved in all inner iterations are achieved by preconditioned CG method
(either with an incomplete Cholesky factorization as an precondtioner or with the Hiptmair-Xu solver \cite{Hiptmair2007}),
with the stopping criterion set to be a relative $l_2$-norm error of the residual less than $10^{-8}$,
if not specified otherwise.
The right-hand side of \eqref{eq.7}, denoted by $b$, is set to be a vector with all components
being ones, and the zero vector is used as the initial guess $x^{(0)}$ for all iterations.

 We shall run respectively the CG  and MINRES   \cite[Algorithms 1 \& 2]{Pestana2013} with
 the new preconditioner $\mathcal{P}$ for solving the saddle-point system~\eqref{eq.7}, and the preconditoned MINRES with the block diagonal preconditioner
 $\mathcal{M}_{\eta, 1/\eta}$, and denote these three methods by
 $\mathcal{P}$-CG,  $\mathcal{P}$-MINRES and  $\mathcal{M}_{\eta, 1/\eta}$-MINRES respectively.
In all our tests, we shall take the parameter $\eta=k^2+1$, unless otherwise defined,
and the outer iterations are terminated based on the criterion
 \begin{equation*}
 \Vert b-\mathcal{K}x^{(k)} \Vert_2\leq  10^{-6}   \cdot \Vert b\Vert_2,
 \end{equation*}
 where $x^{(k)}$ is the $k$th iterate.
 Numbers of iteration by these  solvers with different meshes and wave numbers are listed in Tables \ref{ta.1} and  \ref{ta.2}.
With the new preconditioners, we can see that the required numbers of iteration are slightly smaller than
 the one by the preconditioner $\mathcal{M}_{\eta, 1/\eta}$,
 which is consistent with our theoretical prediction in Section\,\ref{preconditioners}.
 We have also observed that the required numbers of iteration are basically
 independent of mesh size.

In Tables \ref{ta.1} and  \ref{ta.2}, we have also listed the smallest eigenvalue of the matrix
$$A_{\eta}=
\begin{pmatrix}
 A+\eta B^TL^{-1}B-k^2M & 0\\
0 & I_m
\end{pmatrix},$$
 denoted by $\lambda_{\min}(A_{\eta})$, to test its definiteness.
 It is easy to see by \eqref{eq.30} that  we can apply CG with our new preconditioner, under the special inner product defined in \eqref{inpro}  if $A_{\eta}$ is symmetric positive definite.
 It is important to note that for smaller $k$, $ A_{\eta}$ is symmetric positive definite, thus
 CG can be used with the new preconditioner $\mathcal{P}$
 instead of MINRES, though the original system ${\cal K}$ is indefinite.
 The  dotted lines separate  signs of
 $\lambda_{\min}(A_\eta)$ :
 $\lambda_{\min}(A_\eta)$ changes to negative  for $k\geq 1.6$ with Meshes G1 through G4,
 and  for $k\geq 1.25$ with Meshes L1 through L4. This indicates that
 the corresponding preconditioned matrices are no longer  positive definite even under the special non-standard inner product.
 Thus  $\mathcal{P}$-MINRES should be used.
 However,  the numerical results
 indicate that CG  still does not fail even when this violation occurs, and actually converge equally stably and fast.
 This shows very good stability and convergence of CG with the new preconditioner $\mathcal{P}$.
As predicted by Theorem~\ref{tm.bound}, one can  see that the bounds shown by the dotted line is independent of the mesh size.
This is an important feature in applications as it can help us determine which iterative method to use.

\bs
\begin{table}[!htb]\footnotesize\addtolength{\tabcolsep}{-4pt}
	\center
	\caption{Iteration  numbers and values of $\lambda_{\min}(A_{\eta})$ with different $k$ and $\eta=k^2+1$.}
	\label{ta.1}
	\vspace{7pt}
	\begin{tabular}{lccccccc}  
		\hline
		$k$&     0      &    1.0             &{1.55}           &        1.6         &      2  &  4 \\ \hline  
		 {\textbf{Mesh G1}}&  \\
		$\mathcal{P}$-CG&5&	6&	 {11}&	11&	11	&25    \\
		$\mathcal{P}$-MINRES  & 5&	6&	 {11}&	11&	11	&25       \\
		$\mathcal{M}_{\eta,\frac{1}{\eta}}$-MINRES & 6&	9&	 {14}&	13&	13&	30\\
		$\lambda_{\min}(A_{\eta})$ & 0.4677 &   0.4677  &   {0.0340} &  -0.0544 &  -0.8719 &  -7.7031 \\

		 {\textbf{Mesh G2}}&   \\
		$\mathcal{P}$-CG &  5&	7&	 {12}	&12&	11&	25\\
		$\mathcal{P}$-MINRES  &5&	6&	 {12}	&12	&11&	25 \\
		$\mathcal{M}_{\eta,\frac{1}{\eta}}$-MINRES &  6& 	9& 	 {15}& 	15& 	13& 	30    \\
		$\lambda_{\min}(A_{\eta})$ &0.4738  &  0.4738   &  {0.0360}&   -0.0543 &  -0.8988   &-7.9434 \\

		 {\textbf{Mesh G3}}&  \\
		$\mathcal{P}$-CG & 5&	6	& {11}&	11&	11&	25\\
		$\mathcal{P}$-MINRES  &5&	6&	 {11}&	11&	11	&25\\
		$\mathcal{M}_{\eta,\frac{1}{\eta}}$-MINRES &6&	9&	 {15}&	15&	13	&30  \\
		$\lambda_{\min}(A_{\eta})$ &0.4776 &   0.4776&    {0.0369} &  -0.0536 &  -0.8875   &-7.8425 \\

		 {\textbf{Mesh G4}}&  \\
		$\mathcal{P}$-CG &5&	6&	 {9}&	9&	11&	23\\
		$\mathcal{P}$-MINRES  &5&	6&	 {9}&	9&	11&	23 \\
		$\mathcal{M}_{\eta,\frac{1}{\eta}}$-MINRES & 6&	9&	 {11}&	11&	13&	28   \\
		$\lambda_{\min}(A_{\eta})$ &0.4769  &  0.4769  &  {0.0373}&   -0.0533  & -0.8823 &  -7.7907 \\

		 {\textbf{Mesh G5}}&   \\
		$\mathcal{P}$-CG &5& 6 &9 &9 & 11& 23\\
		$\mathcal{P}$-MINRES  &5& 6&9 &9 & 11 &21 \\
		$\mathcal{M}_{\eta,\frac{1}{\eta}}$-MINRES &6 &8 & 11& 11&  13 &28 \\

		\hline
	\end{tabular}
\end{table}

\begin{table}[!htb]\footnotesize\addtolength{\tabcolsep}{-4pt}
	\center
	\caption{Numbers of iterations
		and values of $\lambda_{\min}(A_{\eta})$  with different $k$  and $\eta=k^2+1$.}
	\label{ta.2}
	\vspace{7pt}
	\begin{tabular}{lcccccc}  
		\hline
		$k$&     0   &   1.0       &   {1.2}   & 1.25   &   2       &      4                               \\ \hline  
		{\textbf{Mesh L1}}&   \\
		
		$\mathcal{P}$-CG& 5&	7&  {9}&	8&	12&	25\\
		$\mathcal{P}$-MINRES & 5&	7	& {9}&	8&	10&	24\\
		$\mathcal{M}_{\eta,\frac{1}{\eta}}$-MINRES & 7	&9&	  {11}&	10&	12&	29\\
		
		$\lambda_{\min}(A_{\eta})$ &0.4787  &  0.2496&    { 0.0039} &  -0.0646   &-1.4349 &  -8.3127\\

		{\textbf{Mesh L2}}&  \\
		
		$\mathcal{P}$-CG& 6&	7&	 {9}&	8&	11&	24\\
		$\mathcal{P}$-MINRES & 5&	7& {9}	&8&	11&	24\\
		$\mathcal{M}_{\eta,\frac{1}{\eta}}$-MINRES & 7	&9&	 {11}&	10&	13&	29\\
		$\lambda_{\min}(A_{\eta})$ &0.4582 &   0.2575&   {0.0128}  & -0.0556 &  -1.4249 &  -8.3664\\

		{\textbf{Mesh L3}}&   \\
		
		$\mathcal{P}$-CG&5	&7&	 {9}&	8&	12&	25\\
		$\mathcal{P}$-MINRES &5&	7& {9}	&8&	11&	24\\
		$\mathcal{M}_{\eta,\frac{1}{\eta}}$-MINRES &7&	9	& {11}&	11&	13&	29\\
		$\lambda_{\min}(A_{\eta})$ & 0.4758 &   0.2704  &  {0.0175}&   -0.0530  & -1.4580 &  -8.3974\\

		{\textbf{Mesh L4}}&   \\
		
		$\mathcal{P}$-CG&5 &7&  {8}&8 &10 &24\\
		$\mathcal{P}$-MINRES &5& 7& {8}&8& 10& 24\\
		$\mathcal{M}_{\eta,\frac{1}{\eta}}$-MINRES &6& 9&  {11} &11 &13& 29\\
		$\lambda_{\min}(A_{\eta})$ &0.4654&    0.2753&  {0.0200}  & -0.0512 &  -1.4674 &  -8.4450\\

		{\textbf{Mesh L5}}&   \\
		
		$\mathcal{P}$-CG&5 &7 &8 &8 &10 &24\\
		$\mathcal{P}$-MINRES &5& 7 &8 &8&10 &24\\
		$\mathcal{M}_{\eta,\frac{1}{\eta}}$-MINRES &6 &9& 10& 10& 12&29\\\cline{1-7}
	\end{tabular}
	
\end{table}


In Tables~\ref{ta.3} and \ref{ta.new3}, we test different values  of $\eta-k^2$ to see the influence of the inexact inner solvers on
the required numbers of iteration and to compare the stabilities of the new and
existing preconditioners ${\cal P}$ and $\mathcal{M}_{\eta, 1/\eta}$ when they are
respectively used with CG and MINRES iterations.
For this, we choose a less accurate tolerance $10^{-2}$ for the inner iterations
associated with $L$ and $A+(\eta-k^2)M$.
The required numbers of iteration are reported in Tables~\ref{ta.3} and \ref{ta.new3}
when CG is used with the new and existing preconditioners ${\cal P}$ and
$\mathcal{M}_{\eta, 1/\eta}$ respectively.
As it is expected, the results (lying on the left-hand sides of the dotted lines)
indicate that $\eta-k^2$ should not be too small relatively to $k$
otherwise the matrix $A+(\eta-k^2)M$ becomes nearly singular.
If we ignore those results corresponding to the small values of $\eta-k^2$ relatively to $k$,
we may observe from Tables~\ref{ta.3} and \ref{ta.new3} that
the difference between the required numbers of iterations for these two solvers
increases with $\eta-k^2$.  This indicates better performance and stability of the new preconditioner ${\cal P}$
than the existing one $\mathcal{M}_{\eta, 1/\eta}$.

\begin{table}[!htb]\footnotesize\addtolength{\tabcolsep}{-4pt}
	\center
	\caption{Numbers of iteration for $\mathcal{P}$-CG on Grid~G3  with  different $\eta$, $k$ and inner tolerance $10^{-2}$}
	\label{ta.3}
	\vspace{7pt}
	\begin{tabular}{ccccccccc}
		\hline
		$\eta$&   $k^2+10^{-4}$  &$k^2+1$ & $k^2+4$ & $k^2+8$ & $k^2+20$ & $k^2+45$ \\ \hline
		$k=0$             &  { $>200$ } &  6  &    8  &  10 &   14&  22  \\
		$k=1$           &  {$>200$} & 7& 10  & 12 &  18&26   \\
		$k=2$      & {$>200$}  &  15&     17  &  21 &   28 &  37 \\
		$k=4$          &      $>200$ & $>200$  &  {49} &46&53&  64   \\ \hline
	\end{tabular}
\end{table}

\begin{table}[!htb]\footnotesize\addtolength{\tabcolsep}{-4pt}
	\center
	\caption{Numbers of iteration for $\mathcal{M}_{\eta,\frac{1}{\eta}}$-MINRES on Grid~G3  with different $\eta$, $k$ and inner tolerance $10^{-2}$}
	\label{ta.new3}
	\vspace{7pt}
	\begin{tabular}{cccccccc}
		\hline
		$\eta$&   $k^2+10^{-4}$  &$k^2+1$ & $k^2+4$ & $k^2+8$ & $k^2+20$ & $k^2+45$ \\ \hline
		 {$k=0$}  &       8  &  10    &15  &  17 &   25  &  35  \\
		$k=1$          &    {$>200$} &   11  &  16  &  21  &  28  &  39   \\
		$k=2$          &    {$>200$ }&  23 &   28 &   31 &   40  &  66 \\
		$k=4$   & $>200$ &     {57}  &  50   & 51 &   62  &  80 \\ \hline
	\end{tabular}
\end{table}

\begin{figure}[!hbt]\scriptsize
	\vspace{.7cm}
	\caption{Distributions of eigenvalues smaller than 0.3  of the coefficient matrix $\mathcal{K}$ (upper blue part)  and the matrix $A_{\eta}$        (lower red part) on Grid G3 for $k=4$,  and  $\eta=k^2+1$.}
	\label{fi.fail}
	\includegraphics[width=1\textwidth]{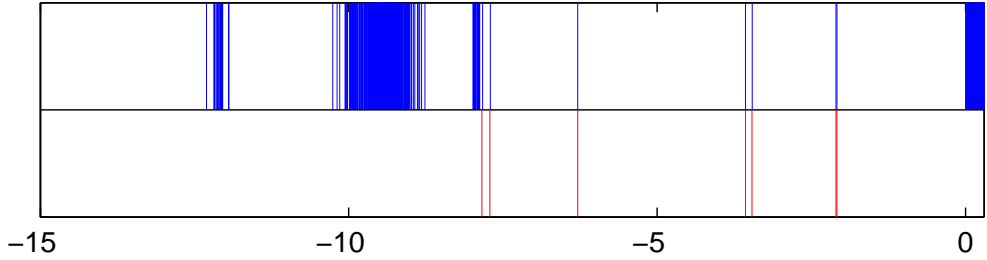}
\end{figure}

We now make some more experiments to further compare the stability of the new and existing
preconditioner $\mathcal{P}$ and $\mathcal{M}_{\eta, 1/\eta}$.
We can clearly see from Table~\ref{ta.3} that CG can be always applied
numerically with the new preconditioner $\mathcal{P}$ and it converges very well,
though it may not guarantee to converge theoretically.
But this is not the case for the preconditioner $\mathcal{M}_{\eta, 1/\eta}$.
%
To see this, we  re-run all the experiments in Table~\ref{ta.new3}, but with CG iteration now
instead MINRES.
In each of the 24 numerical experiments, we have always experienced
the case that one dividend becomes too small, which causes the break-down of the iterative process.
The reasons behind are in fact very simple:
we needs to divide by $p_k^T\mathcal{K}p_k$ (with
$p_k$ being the $k$-th search direction) at the $k$-th CG iteration
with preconditioner $\mathcal{M}_{\eta, 1/\eta}$, and to divide
by $p_k^TA_{\eta} p_k$ at the $k$-th CG iteration
with the new preconditioner ${\cal P}$, due to the existence of a special inner product \eqref{inpro}.
Figure~\ref{fi.fail} shows the distributions of the eigenvalues smaller than $0.3$ of the two matrices
${\cal K}$ and $A_{\eta}$ for $k=4$, and
these smaller and negative eigenvalues contribute mainly to the break-down of the iterations
(most eigenvalues are larger than $0.3$, but not shown in the figure).
As one can see from the table, there are many more eigenvalues in the upper blue part for ${\cal K}$
than in the lower red part for $A_{\eta}$, which explains clearly the highly instability of CG with
the preconditioner $\mathcal{M}_{\eta, 1/\eta}$ and good stability
of CG with the new preconditioner ${\cal P}$.
We have also observed in our experiments that
a larger $\eta-k^2$ makes the eigenvalues $A_{\eta}$ distribute more stably numerically.

\begin{figure}[!hbt]\scriptsize
\vspace{.7cm}
\caption{The eigenvalue distributions of the preconditioned matrix $\mathcal{P}^{-1}\mathcal{K}$ (lower red part)  and  $\mathcal{M}_{\eta,\varepsilon}^{-1}\mathcal{K}$        (upper blue part) on Grid G3 for $k=1.3,$ $\varepsilon=-\frac{1}{\eta-k^2}  $ and  $\eta=k^2+1$.}
\label{fi.c}
\includegraphics[width=1\textwidth]{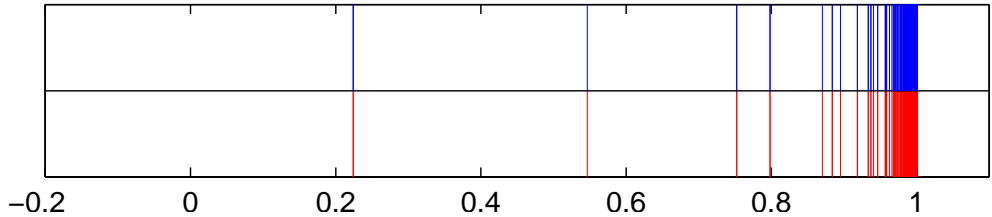}
\end{figure}

\begin{figure}[!thb]

\caption{The eigenvalue distributions of the preconditioned matrix $\mathcal{P}^{-1}\mathcal{K}$ on Grid G3
(from top to bottom: $k=0, 1.3, 1.6, 4$)}
\label{fi.3}

\begin{minipage}[!htb]{1\textwidth}
	\includegraphics[width=1\textwidth]{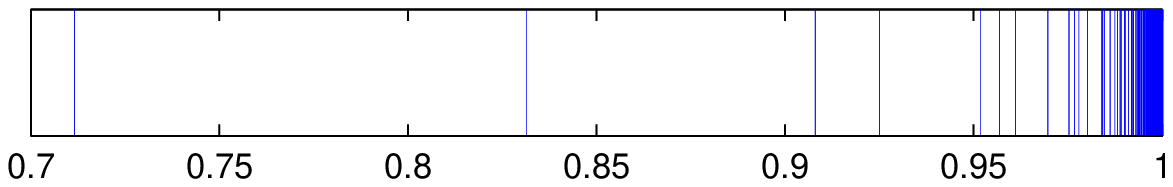}\\
\includegraphics[width=1\textwidth]{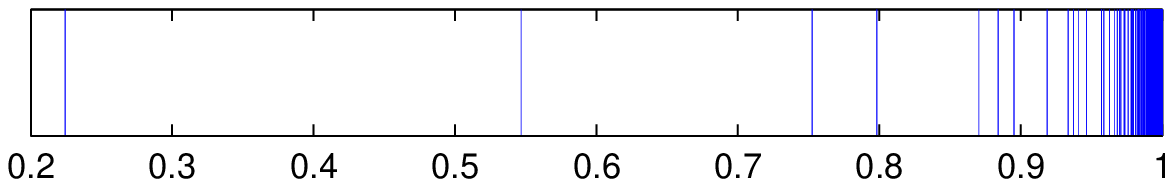}\\
\includegraphics[width=1\textwidth]{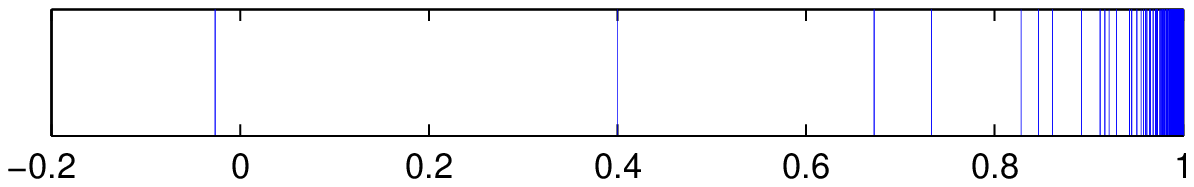}\\
\includegraphics[width=1\textwidth]{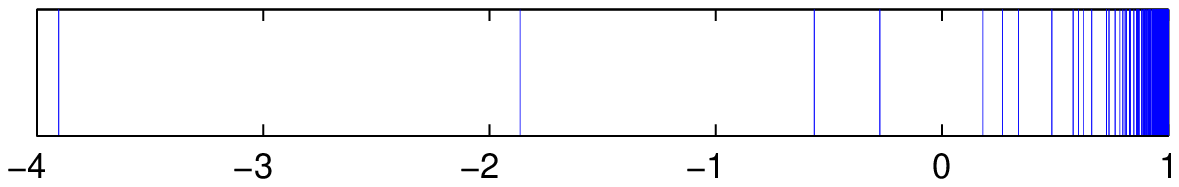}\\
\end{minipage}

\end{figure}

Next we demonstrate the distributions of the eigenvalues of
the preconditioned matrices $\mathcal{P}^{-1}\mathcal{K}$ and $\mathcal{M}_{\eta,\varepsilon}^{-1}\mathcal{K}$.
Figure~\ref{fi.3} plots the eigen-distribution of the preconditioned matrix $\mathcal{P}^{-1}\mathcal{K}$
on Mesh G3 with different wave number $k$.
We can see that
the lower bound of the eigenvalues for $k=1.3$ is about 0.22, and there are $4$ eigenvalues that lie
between 0.22 and $0.8$, while all the remaining 1773 eigenvalues stay in the range $0.8$ and 1, with
850 of them ($m=425$ here) being 1. These results are consistent with our theoretical prediction
(Theorem~\ref{th.6}).
For $k\geq 1.6$,  we see negative eigenvalues, but only a few.
%
Figure~\ref{fi.c} shows that the eigenvalues of  the preconditioned matrices $\mathcal{P}^{-1}\mathcal{K}$   and  $\mathcal{M}_{\eta,\varepsilon}^{-1}\mathcal{K}$   with  $\varepsilon=-\frac{1}{\eta-k^2}$  are exactly the same
for Mesh G3, with $k=1.3$ and  $\eta=k^2+1$.

%

\bs
Finally we test the influence of different right-hand sides $f$ on the required numbers of iteration
when our new preconditioner ${\cal P}$ is used.
Table~\ref{ta.5} shows this influence,
where Df0g, Rf0g and RfRg represent the divergence-free $f$ and $g=0$,
the random $f$ and $g=0$, and the random $f$ and random $g$.
As we can easily see from the table that the  right-hand sides does not cause any effect on
the required numbers of iterations.

\begin{table}[!htb]\footnotesize\addtolength{\tabcolsep}{-0pt}
	\center
	\caption{Numbers of CG iteration on Grid G3 with preconditioner $\mathcal{P}$ and different right-hand sides:
	$k=2$, $\eta=k^2+1$}
				\vspace{4pt}
	\label{ta.5}
		\begin{tabular}{cccccc}
		\hline
		 Grid & G1 & G2 & G3 & G4 & G5 \\ \hline
		Df0g   & 12   & 12   & 12& 12 &12   \\
		Rf0g    & 12   & 12   & 11& 11 &11   \\
		RfRg    & 11   & 12   & 11& 11 &11   \\	\hline
		
	\end{tabular}
\end{table}



\appendix
\section{Appendix}
In this appendix we will generalize the formula for computing the inverse of the symmetric saddle-point matrix
${\cal A}$ in \eqref{eq:inverseA} to the more general case with non-symmtric generalized saddle-point
matrix.
For this purpose, we consider the following non-symmtric generalized saddle-point system:
\begin{equation}	\mathcal{K}\left(\begin{array}{c}u\\ p\end{array}\right) 	\equiv
	\begin{pmatrix}
	A & B \\
	C & D
	\end{pmatrix} \left(\begin{array}{c}u\\ p\end{array}\right)=\left(\begin{array}{c}f\\ g\end{array}\right),
	\end{equation}
where all block matrices $A$, $B$, $C$ and $D$ are allowed to be non-square, with
$A\in \mathbb{R}^{m\times n},$  $B\in \mathbb{R}^{m\times k},$  $C\in \mathbb{R}^{l\times n},$  $D\in \mathbb{R}^{l\times k}$.
But the entire matrix $\mathcal{K}$ is square, i.e., $m+l=n+k=t$.
	
We further assume that the rank of $A$ is $m+n-t$.  Then we write $C_r\in \mathbb{R}^{n\times l}$ as
the matrix of full column rank whose columns span $ker(A)$, and $C_l\in \mathbb{R}^{k\times m}$ as
the matrix of full row rank whose rows span the left null space of $A,$ namely
$C_lA=0$, and $AC_r=0$.
We shall write $L_l=C_lB$ and $L_r=CC_r,$ and
the range space of $A$ by $\mathcal{R}(A)$.
	
We are now ready to formulate the main results in this appendix.
\begin{thm}
	\label{app.1}
Assume that
%
\begin{align}
\mathcal{R}(A)\cap\mathcal{R}(B)={0}\,, & \q \mathcal{R}(A^T) \cap\mathcal{R}(C^T)={0}\,,\label{cond.a}\\
rank(A)=m+n-t\,,& \q rank(B)=k\,, \q rank(C)=l\,,    \label{cond.b}
\end{align}
then the matrix $\mathcal{K}$ is non-singular and its inverse can be represented by
\begin{equation}
\label{appth.1}
\mathcal{K}^{-1}= \begin{pmatrix}
N  & C_rL_r^{-1} \\
L_l^{-1}C_l & 0
\end{pmatrix} ,
\end{equation}
where $N$ satisfies
\begin{align}
NA=I-C_rL_r^{-1}C\,,\q& AN=I-BL_l^{-1}C_l\,,\\
NB=-C_rL_r^{-1} D\,,\q&  CN=-DL_l^{-1}C_l.
\end{align}
If $m=n$, it holds for any $X\in \mathbb{R}^{m\times l}$ such that $A+XC $ is non-singular that
\begin{equation}
N=(A+XC)^{-1}(I-BL_l^{-1}C_l-XDL_l^{-1}C_l).
\end{equation}
\end{thm}

We can easily see that the coefficient matrix in \eqref{eq.2} satisfies   \eqref{cond.a} and \eqref{cond.b}.

For our proof, we first introduce some notations and auxiliary results.
We set $E_A=I-AA^\dagger$ and  $ F_A=I-A^\dagger A$,
where $ A^\dagger$  is the Moore-Penrose inverse of $A,$ namely, $A^\dagger$ is the unique solution $X$ satisfying
 \begin{align}\label{mp}
 XAX=X,\q
 AXA=A,\q
 (AX)^T=AX,\q
 (XA)^T=XA.
 \end{align}
We borrow the following results to work out the explicit formula for the inverse $\mathcal{K}^{-1}$
in  \eqref{appth.1}.
\begin{thm}\cite[Corollary 3.5]{inve}
	\label{app.2}
Under the assumptions of Theorem~\ref{app.1},  $\mathcal{K}$ is non-singular and its inverse is given by
\begin{equation}
\label{appth.2}
\mathcal{K}^{-1}= \begin{pmatrix}
A^\dagger-A^\dagger BB_0^\dagger-C_0^\dagger CA^\dagger-C_0^\dagger\left( D-CA^\dagger B\right)B_0^\dagger  & C_0^\dagger \\
B_0^\dagger & 0
\end{pmatrix} ,
\end{equation}
where $B_0=E_AB$  and $C_0=CF_A$.
\end{thm}

Next we show that the formula \eqref{appth.1} is actually the explicit form of \eqref{appth.2}.
%
%
For this, we first present two auxiliary lemmas.
\begin{lem}
	\label{applem.1}
	We have
	\begin{equation}
	E_A=I-AA^\dagger=C_l^\dagger C_l\,, \q
	F_A=I-A^\dagger A=C_rC_r^\dagger.
	\end{equation}
\end{lem}
{\bf Proof.}
	We prove only the first statement as the second follows similarly.
	Using the Moore-Penrose properties \eqref{mp} it is direct to verify the symmetry of
	$I-AA^\dagger -C_l^\dagger C_l \in\mathbb{R}^{m\times m}$, and
	$$(I-AA^\dagger -C_l^\dagger C_l)A=0,  \q C_l(I-AA^\dagger -C_l^\dagger C_l)=0\,, $$
	which imply that $(I-AA^\dagger -C_l^\dagger C_l)C_l^T=0 $.
	As the  columns of $A$ and $C_l^T$ span $\mathbb{R}^m$,  we immediately see that
	$(I-AA^\dagger -C_l^\dagger C_l)=0 $.

\begin{lem}
		\label{applem.2}
		Under the assumptions of Theorem~\ref{app.1}, the matrices
	$L_l$   and  $L_r$ are non-singular.
\end{lem}
{\bf Proof.}
	We show only the first statement as the second follows similarly.
	Suppose $x^TL_l=x^TC_lB=0$ for $x\in \mathbb{R}^k$, then by direct computing it is easy to see
	$\left(xC_l \,\ 0\right)\mathcal{K}=0,$ so we know $x^TC_l=0$ as ${\cal K}$ is non-singular.
	This implies $x=0$ by noting that $C_l $ is of full row rank. Therefore we know that the square matrix
	$L_l=C_lB$ is non-singular.

Using the above results we can easily derive formulas for the (1,2) and (2,1) blocks  of \eqref{appth.1}.
\begin{thm}\label{lingshi} It holds that
	\begin{align}
	B_0^\dagger  =L_l^{-1}C_l,  \q	C_0^\dagger=C_rL_r^{-1}.
	\end{align}
\end{thm}
{\bf Proof.}
	Again we prove only the first statement.	
	By definitions in Theorem~\ref{app.2} and Lemma~\ref{applem.1} we know $B_0=C_l^\dagger C_lB=C_l^\dagger L_l$.
As the matrix   $L_l^{-1}$ is of full column rank and $C_l$ is of full row rank,
we can derive that $(L_l^{-1}C_l)^\dagger=C_l^\dagger {(L_l^{-1})}^\dagger=B_0$, which ends the desired proof.

Now we go further to derive the explicit formula for computing the (1,1) block of \eqref{appth.2}. For this, we first
show the following results.
\begin{prop} \label{pr.3}
Let $V=A^\dagger-A^\dagger BB_0^\dagger-C_0^\dagger CA^\dagger,$
	$T=V+C_0^\dagger CA^\dagger BB_0^\dagger$, and $N=T-C_0^\dagger DB_0^\dagger$,
	then it holds that
\begin{eqnarray}
&& NA=TA=VA=I-C_rL_r^{-1}C, \q  AN=AT=AV=I-BL_l^{-1}C_l, \label{th1.1}\\
&& TB=CT=0, \q NB=-C_0^\dagger D, \q CN=-DB_0^\dagger. \label{th1.3}
\end{eqnarray}
\end{prop}

{\bf Proof.}
	The first  assertion in \eqref{th1.1} comes from the relations
	\begin{align*}
	VA&=A^\dagger A-A^\dagger BB_0^\dagger A-C_0^\dagger CA^\dagger A
	=(I-C_0^\dagger C)A^\dagger A
	=(I-C_rL_r^{-1} C)(I-C_rC_r^\dagger) \\
	&=I-C_rL_r^{-1} C-C_rC_r^\dagger+C_rL_r^{-1}CC_rC_r^\dagger
	=I-C_rL_r^{-1}C
	\end{align*}
	and  $(C_0^\dagger CA^\dagger BB_0^\dagger)A=0$, which is a direct consequence of the fact that
	$B_0^\dagger A=0$ by Theorem~\ref{lingshi}.
	The second assertion in \eqref{th1.1} follows similarly.
	The relations in \eqref{th1.3} can be readily verified using
	the facts that $B_0^\dagger B=I$ and $CC_0^\dagger=I$, which come from Theorem~\ref{lingshi}.

We end the proof of Theorem~\ref{app.1} with the help of the following results,
which can be checked directly by Proposition~\ref{pr.3}.
\begin{thm}
For any   $X\in \mathbb{R}^{m\times l}$ and $Y\in \mathbb{R}^{k\times n}$, the following identities hold
	\begin{align*}
	N(A+BY)&=I-C_rL_r^{-1}C-C_rL_r^{-1}DY\,,\\
	(A+XC)N&=I-BL_l^{-1}C_l-XDL_l^{-1}C_l\,.
	\end{align*}
\end{thm}

\bibliographystyle{siam}

\end{document}